\documentclass[a4paper,12pt,twocolumn,final]{autart}
\usepackage{color}
\usepackage{lineno}
\usepackage[cmex10]{amsmath}
\usepackage{amscd,amsfonts,amsmath,amssymb,mathrsfs,pifont,stmaryrd,tipa}
\usepackage{array,cases,dsfont,graphicx,texdraw}
\usepackage{graphics} % for pdf, bitmapped graphics files
\usepackage{epsfig} % for postscript graphics files
\usepackage{bm}

\usepackage{geometry}
\usepackage{stfloats}
\usepackage{subfig}
\usepackage{hyperref}

\newtheorem{Remark}{Remark}
\newtheorem{Corollary}{Corollary}
\newtheorem{Definition}{Definition}

\newenvironment{Proof}{\noindent{\em Proof:\/}}{\hfill $\Box$\par}
\newtheorem{Theorem}{Theorem}
\newtheorem{Lemma}{Lemma}

\newtheorem{Assumption}{Assumption}

\begin{document}
\begin{frontmatter}
\title{Distributed Nash Equilibrium Seeking for Games in Systems with Unknown Control Directions}
%: Part \uppercase\expandafter{\romannumeral1} --- Algorithm and Stability Analysis
\thanks[footnoteinfo]{This work was supported by the National Natural Science Foundation of China (NSFC), No. 61803202, the Natural Science Foundation of Jiangsu Province, No. BK20180455, and the Fundamental Research Funds for the Central Universities, No. 30920032203.}
\vspace{-10mm}
\author{Maojiao Ye, Shengyuan Xu, Jizhao Yin}
\thanks{M. Ye, S. Xu and J. Yin are with the School of Automation, Nanjing University of Science and Technology, Nanjing 210094, P.R. China (Email: mjye@njust.edu.cn,syxu@njust.edu.cn,yinjizhao@njust.edu.cn).}


\begin{abstract}
Distributed Nash equilibrium seeking for games in uncertain networked systems without a prior knowledge about control directions is explored in this paper. More specifically, the dynamics of the players are supposed to be first-order or second-order systems in which the control directions are unknown and there are parametric uncertainties. To achieve  Nash equilibrium seeking in a distributed way, Nussbaum function based strategies are proposed through separately designing an optimization module and a state regulation module. The optimization module generates a reference trajectory, that can search for the Nash equilibrium, for the state regulation module. The state regulator is designed to steer the players' actions to the reference trajectory. An adaptive law is included in the state regulation module to compensate for the  uncertain parameter in the players' dynamics and the Nussbaum function is included to address the unavailability of the control directions. Fully distributed implementations of the proposed algorithms are discussed and investigated. Through our analytical explorations, we show that the proposed seeking strategies can drive the players' actions to the Nash equilibrium asymptotically without requiring the homogeneity of the players' unknown control directions based on Barbalat's lemma.  A numerical example is given to support the theoretical analysis of the proposed algorithms.
\end{abstract}
\end{frontmatter}

\section{Introduction}

Games under distributed communication networks  are receiving increasing attention due to their wide applications in numerous  fields. For example, the connectivity control of mobile sensor networks was modeled as a game in which each sensor's objective function contains a local cost that models the sensor's local goal (e.g., source seeking) and a global cost that describes the sensor's willingness to keep connectivity with other sensors \cite{YEcyber2020}. Inspired by the observation that practical engineering systems are usually afflicted with model uncertainties and disturbance, an extended-state-observer based robust Nash equilibrium seeking strategy was proposed in \cite{YEcyber2020}. Energy consumption control might be formulated as an aggregative game, in which each user's cost function depends on the user's own energy consumption and the total energy consumption of all users \cite{Yecyber17}. Dynamic average consensus algorithms can be adapted as an aggregation estimator, based on which distributed Nash equilibrium seeking algorithms were constructed \cite{Yecyber17}. Congestion control problems in wireless sensor networks can be viewed as a semi-aggregative game, in which each data transmitter makes decisions on its data transmission to maximize its own profit \cite{YETAC19}. Interference graphs can be introduced for the interaction descriptions among the data transmitters, based on which Nash equilibrium seeking algorithms were designed in \cite{YETAC19}. Moreover, noncooperative games can be utilized to illustrate the interactions among groups of discrete-time and continuous-time agents under distributed communication networks \cite{MAIJRNC19}. Motivated by the broad applications of networked games, distributed Nash equilibrium seeking has attracted a lot of interests in the past few years and quite a few distributed schemes have been proposed to achieve distributed Nash equilibrium seeking.
%For instance, consensus-based algorithms, extremum seeking based approaches, gossip-based algorithms are typical Nash equilibrium seeking algorithms.
The existing works provide some interesting viewpoints to cope with Nash equilibrium seeking for games in which the players' actions can be freely designed (e.g., \cite{KoshalOR16}\cite{SalehisadaghianiAT16}\cite{YEAT2020}\cite{YEAT18}) or governed by simple dynamics (see e.g., \cite{Ibrahim18}\cite{YETACarXiv}) or possibly subject to disturbance and un-modeled dynamics (see e.g, \cite{YEcyber2020}). A common premise of the existing works is that the control directions are known to the players.

It should be noted that control directions determine the motion directions of a control system and are greatly important as a control force with incorrect direction may deteriorate the system and cause undesired system control performance \cite{ChenAT19}. With the information on control directions, the controller design becomes much simpler. Nevertheless, in some practical circumstances, the control directions are unknown. For instance, due to the inaccurate camera parameters and image depth, the manipulator trajectory tracking control  of visual servo system may need to address the unknown control directions \cite{JiangRobo02}.  Affected by speed variations and loading conditions of the complex, varying environment, the model of ships contains large uncertainties and hence, the autopilot design of time-varying ships requires the accommodation of unknown control directions \cite{DuJOE}. It was recognized that the longitudinal dynamics of the air-breathing hypersonic vehicle suffer from unknown control directions  as well \cite{BunJFI16}. Furthermore, the authors in \cite{Wang2020} argued that in some situations, it is difficult to detect the control directions of quadrotor unmanned aerial vehicles.  Without the information on control directions, the controller design becomes much more challenging especially for multi-agent systems.

Many researchers have been dedicated to investigate systems with unknown control directions. Adaptive designs with Nussbaum-type functions, which can be traced back to \cite{Nussbaum83},  are shown to be effective to deal with uncertainties in control directions.  In \cite{YETAC98}, the Nussbaum-type functions were adopted to achieve adaptive control of nonlinear systems with arbitrary dynamic order and parametric uncertainties. Extremum seekers with unknown control directions were proposed in \cite{ScheinkerTAC13}. Output feedback control for discrete-time systems without a prior control direction knowledge was studied in \cite{YangAT08} in which a discrete Nussbaum gain was utilized to achieve asymptotic output tracking. Nussbaum functions were discussed in \cite{ChenAT19} for systems with time-varying unknown control directions. With the development of multi-agent systems, cooperative control of multi-agent systems with unknown control directions has received increasing attention. For example, the authors in \cite{PengSCL14} considered consensus among a network of first-order integrator-type agents with unknown control directions.  In \cite{ChenTAC17}, the authors supposed that some control directions are known based on which consensus of multi-agent systems with partially unknown and non-identical control directions was addressed. Cooperative output consensus in heterogeneous multi-agent systems with  non-identical control directions was considered in  \cite{GuoTAC17}, where Nussbaum-type functions were adopted to achieve global cooperative output regulation. Distributed optimization among a network of high-order integrator-type agents was addressed in \cite{Tang} without utilizing prior knowledge about control directions. Fully distributed consensus among high-order nonlinear systems in which the agents have heterogenous unknown control directions was investigated in \cite{HuangCyber18}. A new Nussbaum function was employed to deal with the unknown control directions and it was shown that the agents' output can achieve asymptotic consensus \cite{HuangCyber18}.  \textbf{Nevertheless, to the best of the authors' knowledge, distributed Nash equilibrium seeking for networked games in which the players are subject to unknown control directions and uncertain parameters still remains to be addressed.} Motivated by the above observations, this paper tries to shed some light on distributed Nash equilibrium seeking strategy design without utilizing control direction information.

In comparison with the existing works, the main contributions of this paper are summarized as follows.
\begin{enumerate}
  \item Different from the existing works that consider games with known control directions, the seeking strategies proposed in this paper do not require prior direction information. To the best of the authors' knowledge, this is the first work that addresses distributed Nash equilibrium seeking for games with unknown control directions. Besides, this paper also accommodates parametric uncertainties in the players' dynamics.  Through a modular design, this paper proposes Nussbaum function based adaptive seeking strategies to achieve distributed Nash equilibrium seeking for games in both first-order and second-order systems with unknown control directions and parametric uncertainties.
  \item Based on Barbalat's lemma, it is theoretically shown that the players' actions can be steered to the Nash equilibrium while the other auxiliary variables stay bounded by utilizing the proposed algorithms.
  \item Discussions on fully distributed implementation of the proposed algorithms are provided. The explorations show that through adaptive parameter designs, the proposed fully distributed algorithms are effective.
\end{enumerate}

We organize the remaining sections as follows. Some preliminaries are given in Section \ref{ES_NP} and the considered problem is formulated in Section \ref{p1_res}. Section \ref{ne_se_no} presents the main results of the paper, in which first-order and second-order systems with unknown control directions and parametric uncertainties are visited, successively. Discussions on fully distributed implementations of the proposed methods are provided in Section \ref{dis}. Following the theoretical investigations of the developed methods, Section \ref{P_1_exa} provides numerical studies. In the end, conclusions are given in Section \ref{P_1_CONC}.

\section{Preliminaries}\label{ES_NP}

The following definitions or lemmas will be utilized in the rest of the paper.

\begin{Definition}
\cite{PengSCL14} A continuously differentiable function $N_0(\cdot)$ is called a Nussbaum function if
\begin{equation}
\begin{aligned}
&\lim_{q\rightarrow\infty}\sup \frac{1}{q}\int_{0}^{q}N_0(s)ds=\infty,\\
&\lim_{q\rightarrow\infty} \inf \frac{1}{q}\int_{0}^{q}N_0(s)ds=-\infty.
\end{aligned}
\end{equation}
\end{Definition}
Typical examples of Nussbaum functions include $k^2\cos(k)$, $k^2\sin(k)$, to mention just a few. Interested readers are referred to \cite{ChenAT19} for more detailed discussions of Nussbaum functions. In this paper, we adopt $N_0(k)=k^2\sin(k).$
\begin{Lemma}\label{lemma1}
\cite{YETAC98} Suppose that $V(\cdot)$ and $k(\cdot)$ are smooth functions defined on $[0,t_f),$ where $t_f$ is a positive constant and $V(t)\geq0,\forall t\in[0,t_f)$. Moreover, if
\begin{equation}
V(t)\leq \int_0^{t} (a_0 N_0(k(\tau))+1)\dot{k}(\tau)d\tau+c,\forall t\in[0,t_f),
\end{equation}
where $a_0$ is a nonzero constant, $N_0$ is an even smooth Nussbaum function, and $c$ is a suitable constant. Then, $V(t),k(t)$ and $\int_0^{t} (a N_0(k(\tau))+1)\dot{k}(\tau)d\tau$ are bounded on $[0,t_f).$
\end{Lemma}

\begin{Lemma}\label{lemma2}
(\textbf{Barbalat's Lemma} \cite{Slotine}) Suppose that $g(t):\mathbb{R}\rightarrow \mathbb{R}$ is a uniformly continuous function. Then, $\lim_{t\rightarrow \infty}g(t)=0$ given that $\lim_{t\rightarrow \infty} \int_{0}^{t}g(s)ds$ exists and is finite.
\end{Lemma}

A graph $\mathcal{G}$ contains a node set $\mathcal{V}=\{1,2,\cdots,M\}$ ($M\geq 2$ is an integer) and an edge set $\mathcal{E}_d$. The elements of $\mathcal{E}_d$ are represented by $(i,j)$, which illustrates an edge from node $i$ to node $j$ and indicates that node $j$ can receive information from node $i$ but not necessarily vice versa. If $(i,j)\in \mathcal{E}_d$ implies that $(j,i)\in\mathcal{E}_d$ for all $i,j\in\mathcal{V}$. The network is undirected. A directed path from node $i_k$ to node $i_{k+l}$ is a sequence of ordered edges denoted by $(i_{k+j},i_{k+j+1}),j=0,1,2,\cdots,l-1.$
A directed graph is said to be strongly connected if there is a directed path between any two distinct nodes. Similarly, an undirected graph is connected if there is a path between any two distinct nodes. The adjacency matrix $\mathcal{A}$ of a directed graph $\mathcal{G}$ is a matrix whose $(i,j)$th entry is
$a_{ij},$ which is positive if $(j,i)\in \mathcal{E}_d,$ else, $a_{ij}=0.$ Moreover, $a_{ii}=0.$ The adjacency matrix of an undirected graph is similarly defined with a further requirement that $a_{ij}=a_{ji}$ for all $i\neq j.$  Moreover, the Laplacian matrix of graph $\mathcal{G}$ is $\mathcal{L}=\mathcal{D}-\mathcal{A}, $ in which $\mathcal{D}$ is a diagonal matrix whose $i$th diagonal entry is $\sum_{j=1}^{M}a_{ij}$ \cite{YETAC17}\cite{LiTF14}.

\section{Problem Formulation}\label{p1_res}
Consider a game with $N$ players in which the action and cost function of player $i$ is represented by $x_i\in \mathbb{R}$ and $f_i(\mathbf{x}):\mathbb{R}^N\rightarrow \mathbb{R},$  respectively, where $\mathbf{x}=[x_1,x_2,\cdots, x_N]^T$. Denote the player set as $\mathcal{N}=\{1,2,\cdots,N\}$ and suppose that the players' actions are governed by
\begin{equation}\label{dyna_1}
\dot{x}_i=b_iu_i+\phi_i(x_i)\theta_i, \forall i\in\mathcal{N},
\end{equation}
or
\begin{equation}\label{dyna_2}
\begin{aligned}
\dot{x}_i=&v_i,\\
\dot{v}_i=&b_iu_i+\phi_i(x_i)\theta_i, \forall i\in\mathcal{N}.
\end{aligned}
\end{equation}

Note that in \eqref{dyna_1} and \eqref{dyna_2},  $u_i$ is the control input to be designed and $b_i\neq 0$ is an unknown constant. Moreover, $\phi_i(x_i)$ is a sufficiently smooth known function and $\theta_i$ is an unknown parameter. Moreover, $v_i\in\mathbb{R}$ is a state variable of player $i$.

Furthermore, second-order systems in which player $i$'s action is generated by
\begin{equation}\label{dyna_3}
\begin{aligned}
\dot{x}_i=&b_{i1}v_i+\phi_{i1}(x_i)\theta_{i1}\\
\dot{v}_i=&b_{i2}u_i+\phi_{i2}(x_i,v_i)\theta_{i2},\forall i\in\mathcal{N},
\end{aligned}
\end{equation}
where $\theta_{i1},\theta_{i2}, b_{i1},b_{i2}$ are unknown constants, $\phi_{i1}(x_i)$ and $\phi_{i2}(x_i,v_i)$ are smooth functions, will also be considered. Note that in \eqref{dyna_3}, $b_{i1}$ and $b_{i2}$ are nonzero.

The paper aims to design distributed control strategies $u_i$ for systems in
 \eqref{dyna_1}, \eqref{dyna_2} and \eqref{dyna_3}, successively,
 such that $\lim_{t\rightarrow \infty}||\mathbf{x}(t)-\mathbf{x}^*||=0$ where $\mathbf{x}^*$ is the Nash equilibrium defined as follows.
\begin{Definition}\label{def_1}
An action profile $\mathbf{x}^*=(x_i^*,\mathbf{x}_{-i}^*)$ is a Nash equilibrium if for $i\in\mathcal{N},$
\begin{equation}
f_i(x_i^*,\mathbf{x}_{-i}^*)\leq f_i(x_i,\mathbf{x}_{-i}^*),
\end{equation}
for $x_i\in\mathbb{R}$, where $\mathbf{x}_{-i}=[x_1,x_2,\cdots,x_{i-1},x_{i+1},\cdots,x_N]^T$ \cite{YETAC17}.
\end{Definition}
The rest of the paper is based on the following assumptions, which are widely adopted in related works.
\begin{Assumption}\label{ASS1}
For each $i\in\mathcal{N},$ $f_i(\mathbf{x})$ is sufficiently smooth and $\frac{\partial f_i(\mathbf{x})}{\partial x_i}$ is globally Lipshitz with constant $l_i$.
\end{Assumption}
\begin{Assumption}\label{Assu2}
There exists a positive constant $m$ such that for $\mathbf{x},\mathbf{z}\in \mathbb{R}^N,$
\begin{equation}
(\mathbf{x}-\mathbf{z})^T(\mathcal{P}(\mathbf{x})-\mathcal{P}(\mathbf{z}))\geq m||\mathbf{x}-\mathbf{z}||^2,
\end{equation}
where $\mathcal{P}(\mathbf{x})=\left[\frac{\partial f_1(\mathbf{x})}{\partial x_1},\frac{\partial f_2(\mathbf{x})}{\partial x_2},\cdots,\frac{\partial f_N(\mathbf{x})}{\partial x_N}\right]^T.$
\end{Assumption}
\begin{Assumption}\label{ass_3}
The players are equipped with an undirected and connected communication graph $\mathcal{G}$.
\end{Assumption}

For the systems in \eqref{dyna_1} and \eqref{dyna_2}, the nonlinear term should satisfy the following condition.

\begin{Assumption}\label{Assum_4}
For each $i\in\mathcal{N},$ $\phi_i(x_i)$ and $\frac{\partial \phi_
i(x_i)}{\partial x_i}$ are bounded provided that $x_i$ is bounded.
\end{Assumption}

Moreover, for the system in \eqref{dyna_3}, the nonlinear terms should satisfy the following condition.
\begin{Assumption}\label{Assum_5}
For each $i\in\mathcal{N},$ $\phi_{i1}(x_i)$ and $\frac{\partial \phi_
{i1}(x_i)}{\partial x_i}$ are bounded provided that $x_i$ is bounded. Moreover, $\phi_{i2}(x_i,v_i)$ is bounded if $x_i$ and $v_i$ are bounded.
\end{Assumption}
\begin{Remark}
Different from existing works on distributed Nash equilibrium seeking that consider the control directions to be known, we suppose that the control directions are unknown a prior as $b_i$ (or $b_{i1},b_{i2}$) for all $i\in\mathcal{N}$ are not known. Moreover, the players may have different control directions as we do not enforce $sign(b_i)$ (or $sign(b_{i1}),sign(b_{i2})$) for all $i\in\mathcal{N}$ to be the same. Note that in \eqref{dyna_1} and \eqref{dyna_2},  $\theta_i$ is supposed to be unknown as well, indicating that the players are suffering from parametric uncertainties.
\end{Remark}

\section{Main Results}\label{ne_se_no}

In this section, we will establish distributed Nash equilibrium seeking algorithms for games in which the players' actions are governed by \eqref{dyna_1}, \eqref{dyna_2} and \eqref{dyna_3}, successively. In the following, Nash equilibrium seekers that are able to accommodate the unknown control directions and parametric uncertainties will be proposed, followed by their corresponding convergence analyses.
\subsection{Distributed Nash equilibrium seeking for first-order systems with unknown control directions}
In this section, we consider that the action of player $i$ is governed by
\begin{equation}
\dot{x}_i=b_iu_i+\phi_i(x_i)\theta_i, \forall i\in\mathcal{N}.
\end{equation}
In the following, method development and convergence analysis will be presented.
\subsubsection{Method Development}
To achieve distributed Nash equilibrium seeking for systems with unknown control directions, let
\begin{equation}\label{eq1}
u_i=N_0(k_i)(x_i-y_i+\phi_i(x_i)\hat{\theta}_i),
\end{equation}
where $N_0(k_i)=k_i^2sin(k_i)$ and
\begin{equation}\label{eq2}
\begin{aligned}
\dot{k}_i=&(x_i-y_i)(x_i-y_i+\phi_i(x_i)\hat{\theta}_i),\\
\dot{\hat{\theta}}_i=&\phi_i(x_i)(x_i-y_i).
\end{aligned}
\end{equation}
Moreover, $y_i$ is an auxiliary variable generated by
\begin{equation}\label{eq3}
\dot{y}_i=-\nabla_if_i(\mathbf{z}_i),
\end{equation}
where $\nabla_if_i(\mathbf{z}_i)=\frac{\partial f_i(\mathbf{x})}{\partial x_i}\left.\right|_{\mathbf{x}=\mathbf{z}_i}$, $\mathbf{z}_i=[z_{i1},z_{i2},\cdots,z_{iN}]^T$ and
\begin{equation}\label{eq4}
\dot{z}_{ij}=-\delta_{ij}\left(\sum_{k=1}^Na_{ik}(z_{ij}-z_{kj})+a_{ij}(z_{ij}-y_j)\right),
\end{equation}
in which $\delta_{ij}=\delta \bar{\delta}_{ij}$, $\delta$ is positive constant to be determined and
$\bar{\delta}_{ij}$ is a fixed positive constant.

\begin{figure}[!t]
\centering
%\vspace{-20mm}
%\begin{center}
\scalebox{0.6}{\includegraphics[19,459][362,588]{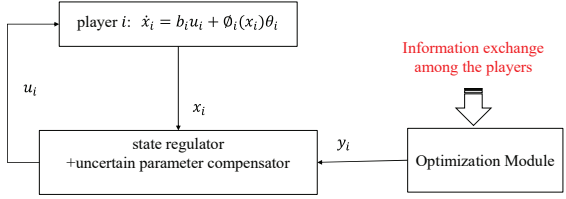}}
%\scalebox{0.35}{\includegraphics{Fig1.eps}}
%\vspace{20mm}
\caption{The illustration of the information flows in the seeking strategy.}\label{fig21}
%\end{center}
\end{figure}

\begin{Remark}
The seeking strategy in \eqref{eq1}-\eqref{eq4} can be viewed as two modules. The subsystem in \eqref{eq1}-\eqref{eq2} is designed to drive $x_i$ to $y_i$. The Nussbaum function in \eqref{eq1} is employed to accommodate the unknown control directions and the second equation in \eqref{eq2} is utilized to compensate the unknown parameter $\theta_i$.
In addition, the subsystem in \eqref{eq3}-\eqref{eq4} is adapted from \cite{YETAC17} to act as  a reference generator that would drive $\mathbf{y}=[y_1,y_2,\cdots,y_N]^T$ to the Nash equilibrium $\mathbf{x}^*$ \cite{YETAC17}.  The schematic outline of \eqref{eq1}-\eqref{eq4} is depicted in the Fig. \ref{fig21}.
\end{Remark}
%\begin{figure}[!htp]
%\centering
%%\vspace{20mm}
%%\begin{center}
%\scalebox{0.36}{\includegraphics[128,226][616,442]{Fig1.eps}}
%%\scalebox{0.35}{\includegraphics{Fig1.eps}}
%%\vspace{55mm}
%\caption{The illustration of the information flows in the seeking strategy.}\label{Fig1}
%%\end{center}
%\end{figure}

\subsubsection{Convergence Analysis}
In this section, we provide the convergence analysis for the seeking strategy proposed in \eqref{eq1}-\eqref{eq4}. Before we proceed to present the convergence results, the following supportive lemma is given.
\begin{Lemma}\label{lemma3}
Suppose that Assumptions \ref{ASS1}-\ref{ass_3} are satisfied. Then, there exists a positive constant $\delta^*$ such that for each $\delta\in(\delta^*,\infty),$ the following conclusions hold:
\begin{itemize}
  \item For each $i,j\in\mathcal{N},$ $y_i(t)$ and $z_{ij}(t)$ are bounded for $t\in[0,\infty).$
  \item For each $i\in\mathcal{N},$  $\dot{y}_i(t)$ globally exponentially decays to zero.
  \item For each $i\in\mathcal{N},$  $\dot{y}_i^2(t)$ is square integrable over $t\in[0,\infty)$, i.e., $\int_{0}^{\infty} \dot{y}_i^2(s)ds\leq c_i$ for some positive constant $c_i.$
\end{itemize}
\end{Lemma}
\begin{Proof}
Following the results in \cite{YETAC17}, it can be obtained that there exists a positive constant $\delta^*$ such that for each $\delta\in(\delta^*,\infty),$  $\mathbf{y}$ and $\mathbf{z}$, where $\mathbf{y}=[y_1,y_2,\cdots,y_N]^T$ and $\mathbf{z}=[\mathbf{z}_1^T,\mathbf{z}_2^T,\cdots,\mathbf{z}_N^T]^T$,  globally exponentially converge to $\mathbf{x}^*$ and $\mathbf{1}_N\otimes \mathbf{x}^*,$ respectively \cite{YETAC17}. Hence, the first conclusion directly follows the results in \cite{YETAC17}. The second conclusion can be reasoned as follows. As $\mathbf{y}$ and $\mathbf{z}$ globally exponentially converge to $\mathbf{x}^*$ and $\mathbf{1}_N\otimes\mathbf{x}^*,$ respectively, there are positive constants $\eta_1,\eta_2$ such that
\begin{equation}
||[(\mathbf{y}-\mathbf{x}^*)^T,(\mathbf{z}-\mathbf{1}_N\otimes \mathbf{x}^*)^T]^T||\leq \eta_1e^{-\eta_2 t}.
\end{equation}

For each $i\in\mathcal{N},$ we get that
\begin{equation}
||\dot{y}_i||=||\nabla_if_i(\mathbf{z}_i)-\nabla_if_i(\mathbf{x}^*)||,
\end{equation}
by noticing that $\nabla_if_i(\mathbf{x}^*)=0,\forall i\in\mathcal{N}$ according to Assumption \ref{Assu2}. By the Lipshitz condition of $\nabla_if_i$ in Assumption \ref{ASS1},  we get that
\begin{equation}
%\begin{aligned}
||\dot{y}_i||\leq l_i||\mathbf{z}_i-\mathbf{x}^*||\leq l_i||\mathbf{z}-\mathbf{1}_N\otimes \mathbf{x}^*||
\leq l_i\eta_1e^{-\eta_2t},
%\end{aligned}
\end{equation}
thus arriving at the second conclusion.

For the third conclusion,
\begin{equation}
\begin{aligned}
\int_{0}^{\infty} \dot{y}_i^2(s)ds=&\int_{0}^{\infty}||\nabla_if_i(\mathbf{z}_i(s))||^2ds\\
\leq & l_i^2 \int_{0}^{\infty}||\mathbf{z}_i(s)-\mathbf{x}^*||^2ds\\
\leq & l_i^2\eta_1^2\int_{0}^{\infty}e^{-2\eta_2s}ds\leq \frac{l_i^2\eta_1^2}{2\eta_2},
\end{aligned}
\end{equation}
thus arriving at the third conclusion with $c_i=\frac{l_i^2\eta_1^2}{2\eta_2}$.
\end{Proof}

Note that by Lemma \ref{lemma3} and \eqref{eq3}-\eqref{eq4}, $\dot{y}_i$ and $\dot{z}_{ij}$ are also bounded as $y_i$ and $z_{ij}$ for all $i,j\in\mathcal{N}$ are bounded.

With the above results in mind, we are now ready to show that the players' actions $\mathbf{x}$ can be driven to the Nash equilibrium $\mathbf{x}^*$ by utilizing the proposed method.

\begin{Theorem}\label{theorem2}
Suppose that Assumptions \ref{ASS1}-\ref{Assum_4} are satisfied. Then, there exists a positive constant $\delta^*$ such that for each $\delta\in(\delta^*,\infty),$
\begin{equation}
\lim_{t\rightarrow\infty}||\mathbf{x}(t)-\mathbf{x}^*||=0,
\end{equation}
and $k_i(t)$, $\hat{\theta}_i(t)$ for all $i\in\mathcal{N}$ stay bounded.
\end{Theorem}
\begin{Proof}
Define a sub-Lyapunov candidate function for player $i$ as
\begin{equation}
V_i=\frac{1}{2}(x_i-y_i)^2+\frac{1}{2}(\theta_i-\hat{\theta}_i)^2.
\end{equation}
Then, the time derivative of $V$ along the trajectory is
\begin{equation}\label{time_de_vi}
\begin{aligned}
\dot{V}_i=&(x_i-y_i)(\dot{x}_i-\dot{y}_i)+(\hat{\theta}_i-\theta_i)\dot{\hat{\theta}}_i\\
=&(x_i-y_i)\left(N_0(k_i)b_i(x_i-y_i+\phi_i(x_i)\hat{\theta}_i)+\phi_i(x_i)\theta_i\right)\\
&-(x_i-y_i)\dot{y}_i+(\hat{\theta}_i-\theta_i)\phi_i(x_i)(x_i-y_i)\\
\leq &N_0(k_i)b_i\dot{k}_i-(x_i-y_i)\dot{y}_i+\hat{\theta}_i\phi_i(x_i)(x_i-y_i)\\
\leq &-(x_i-y_i)^2+(N_0(k_i)b_i+1)\dot{k}_i-(x_i-y_i)\dot{y}_i\\
\leq & -\left(1-\frac{C_i}{2}\right)(x_i-y_i)^2\\
&+(N_0(k_i)b_i+1)\dot{k}_i+\frac{(\dot{y}_i)^2}{2C_i},
\end{aligned}
\end{equation}
by noticing that $|(x_i-y_i)\dot{y}_i|\leq \frac{C_i(x_i-y_i)^2}{2}+\frac{(\dot{y}_i)^2}{2C_i},$ where $C_i$ is a positive constant that satisfies $C_i<2.$

Integrating both sides of \eqref{time_de_vi}, it can be obtained that
 \begin{equation}
 \begin{aligned}
 \int_0^{t_f}\dot{V}_i(\tau)d\tau \leq & -\int_0^{t_f}\left(1-\frac{C_i}{2}\right)(x_i-y_i)^2d\tau\\
 &+\int_0^{t_f}(N_0(k_i)b_i+1)\dot{k}_id\tau+\int_0^{t_f}\frac{(\dot{y}_i)^2}{2C_i}d\tau\\
 \leq &\int_0^{t_f}(N_0(k_i)b_i+1)\dot{k}_id\tau+\frac{c_i}{2C_i}.
 \end{aligned}
 \end{equation}
Note that the last inequality is obtained by noticing that
$\int_0^{t_f}\frac{(\dot{y}_i)^2}{2C_i}d\tau\leq \int_0^{\infty}\frac{(\dot{y}_i)^2}{2C_i}d\tau\leq \frac{c_i}{2C_i}$ according to  Lemma \ref{lemma3}.

Hence,  $V_i(t)$ and $k_i(t)$ are bounded on $[0,t_f)$ by Lemma \ref{lemma1}, which indicates that $x_i-y_i$ and $\hat{\theta}_i$ are bounded. Moreover, as $y_i$ is bounded by Lemma \ref{lemma3}, we obtain that $x_i$ is bounded for $t\in[0,t_f),$ from which we can further obtain that $\dot{x}_i,\dot{k}_i,\dot{\hat{\theta}}_i$ are bounded over the time interval $[0,t_f).$ This implies that there is no finite-time escape for the closed-loop system and hence $t_f=\infty.$

Taking the time derivative of $\dot{k}_i$ gives
\begin{equation}
\begin{aligned}
\ddot{k}_i=&(\dot{x}_i-\dot{y}_i)(x_i-y_i+\phi_i(x_i)\hat{\theta}_i)\\
&+(x_i-y_i)(\dot{x}_i-\dot{y}_i+\frac{\partial \phi_i(x_i)}{\partial x_i}\dot{x}_i\hat{\theta}_i+\phi_i(x_i)\dot{\hat{\theta}}_i).
\end{aligned}
\end{equation}

As $x_i$ is bounded for $t\in[0,\infty)$, $\frac{\partial \phi_i(x_i)}{\partial x_i}$ is bounded  for $t\in[0,\infty)$ by Assumption \ref{Assum_4}.  Moreover, noticing that
$x_i,y_i,\phi_i(x_i),\hat{\theta}_i,\dot{x}_i,\dot{y}_i,\dot{\hat{\theta}}_i$ are all bounded, we get that $\ddot{k}_i$ is bounded. Hence, $\dot{k}_i(t)$ is uniformly continuous with respect to $t.$ In addition,
\begin{equation}
%\begin{aligned}
 \int_0^{\infty}\dot{k}_i(s)ds= k_i(\infty)-k_i(0)\leq k_i^*,
%\end{aligned}
\end{equation}
where $k_i^*$ is a finite constant determined by the bounds of $k_i(t).$

Therefore, $(x_i(t)-y_i(t))(x_i(t)-y_i(t)+\phi_i(x_i(t))\hat{\theta}_i(t))$ is integrable over $t\in[0,\infty)$.  Hence
\begin{equation}\label{end_1}
\lim_{t\rightarrow\infty}(x_i(t)-y_i(t))(x_i(t)-y_i(t)+\phi_i(x_i(t))\hat{\theta}_i(t))=0,
\end{equation}
by Lemma \ref{lemma2}.

From the other aspect, taking the time derivative of $\dot{\hat{\theta}}_i$ gives
\begin{equation}
\ddot{\hat{\theta}}_i=\frac{\partial \phi_i(x_i)}{\partial x_i}\dot{x}_i(x_i-y_i)+\phi_i(x_i)(\dot{x}_i-\dot{y}_i),
\end{equation}
from which we see that $\ddot{\hat{\theta}}_i$ is bounded by noticing that $x_i,y_i,\dot{x}_i,\dot{y}_i,\frac{\partial \phi_i(x_i)}{\partial x_i},\phi_i(x_i)$ are bounded.

Therefore,  $\dot{\hat{\theta}}_i$ is uniformly continuous with respect to $t$. Moreover,
\begin{equation}
\int_0^{\infty}\dot{\hat{\theta}}_i(t)dt=\hat{\theta}_i(\infty)-\hat{\theta}_i(0)\leq \hat{\theta}_i^*,
\end{equation}
where $\hat{\theta}_i^*$ is a constant determined by the bounds of $\hat{\theta}_i.$ Hence, by Lemma \ref{lemma2}, we can obtain that
\begin{equation}\label{end_2}
\lim_{t\rightarrow\infty} \phi_i(x_i(t))(x_i(t)-y_i(t))=0.
\end{equation}
By \eqref{end_1}, we have
$x_i(t)=y_i(t)$ or alternatively, $x_i(t)-y_i(t)+\phi_i(x_i(t))\hat{\theta}_i(t)=0$ for $t=\infty.$ Moreover, by \eqref{end_2}, we have $\phi_i(x_i(t))=0$ or $x_i(t)=y_i(t)$ for $t=\infty.$
Suppose that $x_i(t)\neq y_i(t)$ for $t=\infty,$ then, $\phi_i(x_i(t))=0$ must be satisfied. If this is the case, $x_i(t)-y_i(t)+\phi_i(x_i(t))\hat{\theta}_i(t)=x_i(t)-y_i(t)\neq 0$ for $t=\infty,$ indicating that  \eqref{end_1} can not be satisfied. Hence, we arrive at a contradiction and obtain that $x_i(t)= y_i(t)$ must be satisfied for $t=\infty.$ Recalling that $\mathbf{y}(t)\rightarrow \mathbf{x}^*$ as $t\rightarrow\infty$, which is proven in Lemma \ref{lemma3}, we arrive at the conclusion that $\mathbf{x}(t)\rightarrow \mathbf{x}^*$ as $t\rightarrow\infty$, thus completing the proof.
\end{Proof}

In this paper, we focus on the case in which the communication graph is undirected and connected for simplicity. However, it should be noted that the presented results are still valid for strongly connected digraphs. To highlight this point, the following corollary is given.
\begin{Corollary}\label{coro2}
Suppose that Assumptions \ref{ASS1}-\ref{Assu2}, \ref{Assum_4} are satisfied and the communication graph is strongly connected. Then, there exists a positive constant $\delta^*$ such that for each $\delta\in(\delta^*,\infty),$
\begin{equation}
\lim_{t\rightarrow\infty}||\mathbf{x}(t)-\mathbf{x}^*||=0,
\end{equation}
and $k_i(t),\hat{\theta}_i(t)$ for all $i\in\mathcal{N}$ stay bounded.
\end{Corollary}
\begin{Proof}
The proof follows that of Theorem \ref{theorem2} by noticing that the results in Lemma \ref{lemma3} are still valid for strongly connected digraphs.
\end{Proof}

In Theorem \ref{theorem2}, we consider that each player $i$'s action is subject to both unknown control directions ($b_i$ is unknown) and uncertain parameter $\theta_i,$ i.e.,
\begin{equation}
\dot{x}_i=b_iu_i+\phi_i(x_i)\theta_i.
\end{equation}
If there is no uncertain parameter, and the players' actions are generated by
\begin{equation}
\dot{x}_i=b_iu_i.
\end{equation}
Then, the proposed seeking strategy can be revised to be
\begin{equation}\label{eq11}
u_i=N_0(k_i)(x_i-y_i),
\end{equation}
where $N_0(k_i)=k_i^2sin(k_i)$,
\begin{equation}\label{eq12}
\dot{k}_i=(x_i-y_i)^2,
\end{equation}
and $y_i$ is generated by \eqref{eq3}-\eqref{eq4}. If this is the case, the following corollary can be obtained.

\begin{Corollary}
Suppose that Assumptions \ref{ASS1}-\ref{ass_3} are satisfied. Then, there exists a positive constant $\delta^*$ such that for each $\delta\in(\delta^*,\infty),$
\begin{equation}
\lim_{t\rightarrow\infty}||\mathbf{x}(t)-\mathbf{x}^*||=0,
\end{equation}
and $k_i(t)$ for all $i\in\mathcal{N}$ stay bounded.
\end{Corollary}

\subsection{Distributed Nash equilibrium seeking for second-order systems}
In this section, we suppose that for each $i\in\mathcal{N},$ player $i$'s action $x_i$ is governed by
\begin{equation}\label{seca}
\begin{aligned}
\dot{x}_i=&v_i\\
\dot{v}_i=&b_iu_i+\phi_i(x_i)\theta_i,
\end{aligned}
\end{equation}
in which $v_i\in\mathbb{R}$ is a state of player $i$.
\subsubsection{Method development}
To achieve distributed Nash equilibrium seeking for games in which each player $i$'s dynamics is governed by \eqref{seca}, the control input $u_i$ is designed as
\begin{equation}\label{eqq1}
u_i=N_0(k_i)(x_i-y_i+v_i+\phi_i(x_i)\hat{\theta}_i+(\dot{x}_i-\dot{y}_i)),
\end{equation}
where $N_0(k_i)=k_i^2sin(k_i)$ and
\begin{equation}\label{eqq2}
\begin{aligned}
\dot{k}_i=&(x_i-y_i+v_i)(x_i-y_i+v_i+\phi_i(x_i)\hat{\theta}_i+(\dot{x}_i-\dot{y}_i)),\\
\dot{\hat{\theta}}_i=&\phi_i(x_i)(x_i-y_i+v_i).
\end{aligned}
\end{equation}
Moreover, $y_i$ is an auxiliary variable generated by
\begin{equation}\label{eqq3}
\dot{y}_i=-\nabla_if_i(\mathbf{z}_i),
\end{equation}
where $\mathbf{z}_i=[z_{i1},z_{i2},\cdots,z_{iN}]^T$. Furthermore,
\begin{equation}\label{eqq4}
\dot{z}_{ij}=-\delta_{ij}\left(\sum_{k=1}^Na_{ik}(z_{ij}-z_{kj})+a_{ij}(z_{ij}-y_j)\right),
\end{equation}
where $\delta_{ij}=\delta \bar{\delta}_{ij}$, $\delta$ is positive constant to be determined and
$\bar{\delta}_{ij}$ is a fixed positive constant.

To establish the results for second-order systems, the following assumption is also needed.
\begin{Assumption}\label{ass_5}
For each $i,j\in\mathcal{N},$ $\frac{\partial \nabla_i f_i(\mathbf{x})}{\partial x_j}$ is bounded given that $\mathbf{x}$ is bounded.
\end{Assumption}
\begin{Remark}
Compared the strategy in \eqref{eqq1}-\eqref{eqq4} with \eqref{eq1}-\eqref{eq4}, we see that the optimization modules are the same while the regulation modules are different. As the system in \eqref{seca} is a second-order system, we further utilize $\dot{x}_i$ and $\dot{y}_i$ in the seeking strategy. Recalling the definitions of $\dot{x}_i$ and $\dot{y}_i$, it is clear that the communication in the proposed seeking strategy is still one-hop.
\end{Remark}
\subsubsection{Convergence analysis}

The following theorem illustrates the convergence result for the proposed method.
\begin{Theorem}\label{th2}
Suppose that Assumptions \ref{ASS1}-\ref{ass_3}, \ref{Assum_5}-\ref{ass_5} are satisfied. Then, there exists a positive constant $\delta^*$ such that for each $\delta\in(\delta^*,\infty),$
\begin{equation}
\lim_{t\rightarrow\infty} ||\mathbf{x}(t)-\mathbf{x}^*||\rightarrow 0.
\end{equation}
Moreover, $k_i(t)$ and $\hat{\theta}_i(t)$ for all $i\in\mathcal{N}$ stay bounded.
\end{Theorem}
\begin{Proof}
For notational convenience, let $\xi_i=x_i-y_i+v_i$. Define the sub-Lyapunov candidate function for player $i$ as
\begin{equation}
V_i=\frac{1}{2}(x_i-y_i)^2+\frac{1}{2}\xi_i^2+\frac{1}{2}(\theta_i-\hat{\theta}_i)^2.
\end{equation}
Then, the time derivative of $V_i$ is
\begin{equation}\label{int}
\begin{aligned}
\dot{V}_i=&-(x_i-y_i)(x_i-y_i-\xi_i)-(x_i-y_i)\dot{y}_i\\
&+\xi_i\left(b_iN_0(k_i)(\xi_i+\phi_i(x_i)\hat{\theta}_i+(\dot{x}_i-\dot{y}_i))\right)\\
&+\xi_i\left(\phi_i(x_i)\theta_i+\dot{x}_i-\dot{y}_i\right)+(\hat{\theta}_i-\theta_i)\phi_i(x_i)\xi_i\\
=&-(x_i-y_i)^2-\xi_i^2+(b_iN_0(k_i)+1)\dot{k}_i\\
&+(x_i-y_i)\xi_i-(x_i-y_i)\dot{y}_i\\
\leq &-\left(\frac{1}{2}-\frac{C_i}{2}\right)(x_i-y_i)^2-\frac{1}{2}\xi_i^2\\
&+(b_iN_0(k_i)+1)\dot{k}_i+\frac{(\dot{y}_i)^2}{2C_i},
\end{aligned}
\end{equation}
where $C_i$ is a positive constant that satisfies $C_i<1.$

Integrating both sides of \eqref{int} over $t\in[0,t_f)$ gives
 \begin{equation}
 \begin{aligned}
 \int_0^{t_f}\dot{V}_id\tau \leq & -\int_0^{t_f}\left[\left(\frac{1}{2}-\frac{C_i}{2}\right)(x_i-y_i)^2+\frac{1}{2}\xi_i^2\right]d\tau\\
&+\int_0^{t_f}(b_iN_0(k_i)+1)\dot{k}_id\tau+\int_0^{t_f}\frac{(\dot{y}_i)^2}{2C_i}d\tau\\
\leq & \int_0^{t_f}(b_iN_0(k_i)+1)\dot{k}_id\tau+ \frac{c_i}{2C_i}.
\end{aligned}
\end{equation}

Hence, by Lemma \ref{lemma1}, we get that $V_i$ and $k_i$ are bounded for $t\in[0,t_f),$ which further indicates that $x_i-y_i,\xi_i,\hat{\theta}_i$ are bounded for $t\in[0,t_f)$. Recalling that $y_i$ is bounded, we get that $x_i$ is bounded for $t\in[0,t_f)$. Hence, $v_i$ is bounded. Therefore, there is no finite-time escape for the closed-loop system, which indicates that $t_f=\infty.$ Recalling  \eqref{eqq1}-\eqref{eqq4}, we can obtain that $\dot{x}_i,\dot{v}_i,\dot{y}_i,\dot{k}_i,\dot{\hat{\theta}}_i$ are all bounded.  Taking the time derivative of $\dot{k}_i(t)$ gives
\begin{equation}
\begin{aligned}
\ddot{k}_i=&(\dot{v}_i+\dot{x}_i-\dot{y}_i)(\xi_i+\phi_i(x_i)\hat{\theta}_i+(\dot{x}_i-\dot{y}_i))\\
&+\xi_i(\dot{x}_i-\dot{y}_i+\dot{v}_i+\frac{\partial \phi_i(x_i)}{\partial x_i}\dot{x}_i\hat{\theta}_i+\phi_i(x_i)\dot{\hat{\theta}}_i)\\
&+\xi_i(\ddot{x}_i-\ddot{y}_i).\\
\end{aligned}
\end{equation}
Note that $\ddot{x}_i=\dot{v}_i$ is bounded and $\ddot{y}_i=\left(\frac{\nabla_i f_i(\mathbf{x})}{\partial \mathbf{x}}\left.\right|_{\mathbf{x}=\mathbf{z}_i}\right)^T \dot{\mathbf{z}}_i$ is bounded as $\mathbf{z}_i$, $\dot{\mathbf{z}}_i$ are bounded and $\frac{\nabla_i f_i(\mathbf{x})}{\partial \mathbf{x}}\left.\right|_{\mathbf{x}=\mathbf{z}_i}$ is bounded for bounded $\mathbf{z}_i$ (by Assumption \ref{ass_5}),  it can be seen that $\ddot{k}_i$ is bounded. Hence, $\dot{k}_i(t)$ is uniformly continuous. Moreover,
\begin{equation}
\int_0^{\infty}\dot{k}_i(\tau)d\tau=k_i(\infty)-k_i(0)\leq k_i^*,
\end{equation}
where $k_i^*$ is a constant determined by the bounds of $k_i(t).$
Hence, by Lemma \ref{lemma2}, we can obtain that $(v_i+x_i-y_i)(x_i-y_i+v_i+\phi_i(x_i)\hat{\theta}_i+(\dot{x}_i-\dot{y}_i))\rightarrow 0$ as $t\rightarrow\infty.$

Similarly,
\begin{equation}
\int_0^{\infty}\dot{\hat{\theta}}_i(\tau)d\tau=\hat{\theta}_i(\infty)-\hat{\theta}_i(0)\leq \hat{\theta}_i^*,
\end{equation}
where $\hat{\theta}_i^*$ is a constant determined by the bounds of $\hat{\theta}_i.$ Hence, by Lemma \ref{lemma2},  we can obtain that $\phi_i(x_i)(v_i+x_i-y_i)\rightarrow 0$ as $t\rightarrow \infty.$

Hence, for $t=\infty,$ we have $\phi_i(x_i)(v_i+x_i-y_i)=0$ and $(v_i+x_i-y_i)(x_i-y_i+v_i+\phi_i(x_i)\hat{\theta}_i+(\dot{x}_i-\dot{y}_i))=0.$

Case I: $\phi_i(x_i)=0$ but $v_i+x_i-y_i\neq 0$ for $t=\infty.$
In this case,  $x_i-y_i+v_i+(\dot{x}_i-\dot{y}_i)=0$. Recalling that as $t\rightarrow\infty,$ $y_i\rightarrow x_i^*,$  and $\dot{y}_i\rightarrow 0$, we get that
\begin{equation}
\dot{x}_i=-\frac{1}{2}(x_i-x_i^*),
\end{equation}
from which it is clear that $\mathbf{x}(t)\rightarrow \mathbf{x}^*$ for $t\rightarrow \infty.$

Case II: $v_i+x_i-y_i=0$ for $t=\infty.$ If this is the case
\begin{equation}
\dot{x}_i=-(x_i-y_i),
\end{equation}
as $t\rightarrow \infty.$ Recalling that $\lim_{t\rightarrow\infty} (y_i(t)-x_i^*)=0,$ we can obtain that $\lim_{t\rightarrow\infty}||\mathbf{x}(t)-\mathbf{x}^*||=0.$ To this end, the conclusion is obtained.
\end{Proof}

Similar to Corollary \ref{coro2}, the following result can be obtained if the communication graph is strongly connected.
\begin{Corollary}\label{coro4}
Suppose that Assumptions \ref{ASS1}-\ref{Assu2}, \ref{Assum_4}, \ref{ass_5} are satisfied and the communication graph is strongly connected. Then, there exists a positive constant $\delta^*$ such that for each $\delta\in(\delta^*,\infty),$
\begin{equation}
\lim_{t\rightarrow\infty}||\mathbf{x}(t)-\mathbf{x}^*||=0
\end{equation}
and $k_i(t),\hat{\theta}_i(t)$ for all $i\in\mathcal{N}$ stay bounded.
\end{Corollary}

\subsection{Distributed Nash equilibrium seeking for more general second-order systems}
In this section, we consider a game in which each player $i$'s action is governed by
\begin{equation}\label{thri}
\begin{aligned}
\dot{x}_i=&b_{i1}v_i+\phi_{i1}(x_i)\theta_{i1}\\
\dot{v}_i=&b_{i2}u_i+\phi_{i2}(x_i,v_i)\theta_{i2}.
\end{aligned}
\end{equation}
%where $\theta_{i1},\theta_{i2}, b_{i1},b_{i2}$ are unknown constants, $\phi_{i1}(x_i)$ and $\phi_{i2}(x_i,v_i)$ are smooth functions. Moreover, $b_{i1}$ and $b_{i2}$ are nonzero.
\begin{Remark}
Note that compared with \eqref{seca}, the effect of $v_i$ on $x_i$ is also uncertain in \eqref{thri} as $b_{i1}$ is also unknown. In addition, an uncertain nonlinear term $\phi_{i1}(x_i)\theta_{i1}$ is addressed as well. Hence, in this problem, $b_{i1}$ and $b_{i2}$ are both unknown directions that should be addressed. Moreover, both $\theta_{i1}$ and $\theta_{i2}$ result in uncertain nonlinearities that should be accommodated.
\end{Remark}

Motivated by \cite{YETAC98}, the Nash equilibrium seeking strategy is designed in the following process:

\textbf{Step 1:} Generate a reference trajectory $y_i$ for $i\in\mathcal{N}$ that would converge to the Nash equilibrium according to
 \begin{equation}\label{eqq13}
 \begin{aligned}
\dot{y}_i=&-\nabla_if_i(\mathbf{z}_i)\\
\dot{z}_{ij}=&-\delta_{ij}\left(\sum_{k=1}^Na_{ik}(z_{ij}-z_{kj})+a_{ij}(z_{ij}-y_j)\right),
\end{aligned}
\end{equation}
where $j\in\mathcal{N},$ $\mathbf{z}_i=[z_{i1},z_{i2},\cdots,z_{iN}]^T$, $\delta_{ij}=\delta \bar{\delta}_{ij}$, $\delta$ is positive constant to be determined and
$\bar{\delta}_{ij}$ is a fixed positive constant.

\textbf{Step 2:} Generate a reference trajectory $\alpha_i$ for $v_i$ as
\begin{equation}\label{eqq14}
\begin{aligned}
\alpha_i=&N_0(k_{i1})(x_i-y_i+\phi_{i1}(x_i)\hat{\theta}_{i1})\\
\dot{k}_{i1}=&(x_i-y_i)(x_i-y_i+\phi_{i1}(x_i)\hat{\theta}_{i1})\\
\dot{\hat{\theta}}_{i1}=&\phi_{i1}(x_i)(x_i-y_i).
\end{aligned}
\end{equation}

\textbf{Step 3:} Let $\beta_i=v_i-\alpha_i$. Then, through direct calculation, it can be obtained that
\begin{equation}
\begin{aligned}
\dot{\beta}_i=&\dot{v}_i-\dot{\alpha}_i\\
=&b_{i2}u_i+\phi_{i2}(x_i,v_i)\theta_{i2}+\Psi_{i1}(x_i,k_i,\hat{\theta}_{i1})\theta_{i1}\\
&+\Psi_{i2}(k_{i1},x_i,y_i,\hat{\theta}_{i1},\dot{y}_i)+\Psi_{i3}(k_i,x_i,\hat{\theta}_{i1},v_i)b_{i1}
\end{aligned}
\end{equation}
where $\Psi_{i1}=-N_0(k_{i1})\left(\phi_{i1}(x_i)+\frac{\partial \phi_{i1}(x_i)}{\partial x_i}\hat{\theta}_{i1}\phi_{i1}(x_i)\right),$ $\Psi_{i2}=-(2k_{i1}sin(k_{i1})+k_{i1}^2cos(k_{i1}))(x_i-y_i)(x_i-y_i+\phi_{i1}(x_i)\hat{\theta}_{i1})^2-N_0(k_{i1})(-\dot{y}_i+\phi_{i1}^2(x_i)(x_i-y_i))$ and $\Psi_{i3}=-N_0(k_{i1})\left(\frac{\partial \phi_{i1}(x_i)}{\partial x_i}\hat{\theta}_{i1}+1\right)v_i.$
Accordingly, the control input $u_i$ is designed as
\begin{equation}\label{eqq15}
\begin{aligned}
u_i=&N_0(k_{i2})(\beta_i+\phi_{i2}\bar{\theta}_{i2}+\Psi_{i1}\bar{\theta}_{i1}+\Psi_{i2}+\Psi_{i3}\bar{b}_{i1}),\\
\dot{k}_{i2}=&\beta_i(\beta_i+\phi_{i2}\bar{\theta}_{i2}+\Psi_{i1}\bar{\theta}_{i1}+\Psi_{i2}+\Psi_{i3}\bar{b}_{i1}),\\
\dot{\bar{\theta}}_{i2}=&\beta_i \phi_{i2},\dot{\bar{\theta}}_{i1}=\beta_i \Psi_{i1},\\
\dot{\bar{b}}_{i1}=&\beta_i \Psi_{i3}.
\end{aligned}
\end{equation}

\begin{Remark}
Note that \eqref{eqq14} is designed to drive $x_i$ to $y_i$ and \eqref{eqq15} is designed to drive $v_i$ to $\alpha_i.$ The design of the control input in \eqref{eqq15} is motivated by  \cite{YETAC98} that treats $v_i$ as a virtual control input for $x_i.$ To deal with unknown constants $b_{i1}$ and $b_{i2},$ two Nussbaum functions are included. To accommodate multiple Nussbaum functions, the idea is to design the control input such that $\beta_i$ is square integrable (see also \cite{YETAC98}).
\end{Remark}

The following theorem establishes the stability of Nash equilibrium under the control input designed in \eqref{eqq15}.
\begin{Theorem}
Suppose that Assumptions \ref{ASS1}-\ref{ass_3},\ref{Assum_5} are satisfied. Then, there exists a positive constant $\delta^*$ such that for each $\delta\in(\delta^*,\infty),$
\begin{equation}
\lim_{t\rightarrow\infty}||\mathbf{x}(t)-\mathbf{x}^*||=0,
\end{equation}
 and other variables stay bounded.
\end{Theorem}
\begin{Proof}
The proof is similar to those in \cite{YETAC98} and Theorem \ref{th2}. For the convenience of the readers, sketch of the proof is given as follows.

Step 1: Show that $\beta_i$ is square integrable by defining the sub-Lyapunov candidate function as
\begin{equation}
\begin{aligned}
V_{i1}=&\frac{1}{2}\beta_i^2+\frac{1}{2}(\bar{\theta}_{i2}-\theta_{i2})^2\\
&+\frac{1}{2}(\bar{\theta}_{i1}-\theta_{i1})^2+\frac{1}{2}(\bar{b}_{i1}-b_{i1})^2.
\end{aligned}
\end{equation}
Then, following the proof of Theorem \ref{th2}, it can be obtained that
\begin{equation}\label{int_ff}
\dot{V}_{i1}\leq-\beta_i^2+ (b_{i2}N_0(k_{i2})+1)\dot{k}_{i2}.
\end{equation}

Moreover, taking integrations on both sides of \eqref{int_ff} over $[0,t_f)$, we get that
\begin{equation}\label{int_ff1}
\int_{0}^{t_f}\dot{V}_{i1}d\tau \leq-\int_{0}^{t_f}\beta_i^2d\tau+ \int_{0}^{t_f}(b_{i2}N_0(k_{i2})+1)\dot{k}_{i2}d\tau,
\end{equation}
from which it can be obtained that $V_{i1}$, $k_{i2}$ and $\int_0^{t_f}(b_{i2}N_0(k_{i2})+1)\dot{k}_{i2}dt$ are bounded by Lemma \ref{lemma1}.

Moreover, from \eqref{int_ff1}, it is clear that
\begin{equation}
\int_{0}^{t_f}\beta_i^2d\tau\leq V_i(0)-V_i(t_f)+ \int_{0}^{t_f}(b_{i2}N_0(k_{i2})+1)\dot{k}_{i2}d\tau.
\end{equation}
Hence, $\beta_i$ is square integrable for $t\in[0,t_f)$.

Step 2: Show that $x_i$ can be driven to $y_i$ by defining the other sub-Lyapunov function as
\begin{equation}
V_{i2}=\frac{1}{2}(x_i-y_i)^2+(\hat{\theta}_{i1}-\theta_{i1})^2.
\end{equation}

Then, the time derivative of $V_{i2}$ is
\begin{equation}
\begin{aligned}
\dot{V}_{i2}=&(x_i-y_i)(b_{i1}\alpha_i+b_{i1}\beta_i+\phi_{i1}(x_i)\theta_{i1}-\dot{y}_i)\\
&+(\hat{\theta}_{i1}-\theta_{i1})\dot{\hat{\theta}}_{i1}\\
=&(x_i-y_i)(b_{i1}N_0(k_{i1})(x_i-y_i+\phi_{i1}(x_i)\hat{\theta}_{i1})\\
&+b_{i1}\beta_i+\phi_{i1}(x_i)\theta_{i1}-\dot{y}_i)+(\hat{\theta}_{i1}-\theta_{i1})\dot{\hat{\theta}}_{i1}\\
\leq &-(x_i-y_i)^2+(b_{i1}N_0(k_{i1})+1)\dot{k}_{i1}\\
&+(x_i-y_i)b_{i1}\beta_i-(x_i-y_i)\dot{y}_i\\
\leq & -\frac{1}{2}(x_i-y_i)^2+(b_{i1}N_0(k_{i1})+1)\dot{k}_{i1}\\
&+\frac{C_{i1}}{2}b_{i1}^2\beta_i^2+\frac{C_{i2}}{2}\dot{y}_i^2.
\end{aligned}
\end{equation}
where $C_{i1},C_{i2}$ are positive constants that satisfy $\frac{1}{C_{i1}}+\frac{1}{C_{i2}}\leq 1$. Noticing that both $\beta_i$ and $\dot{y}_i$ are square integrable for $t\in[0,t_f),$ we obtain that $V_{i1},k_{i1}$ are bounded for $t\in[0,t_f)$. Combining the above two steps, it can be seen that $x_i-y_i,$ $k_{i1},$ $\hat{\theta}_{i1}$ as well as $\beta_i$, $k_{i2},\bar{\theta}_{i1},\bar{\theta}_{i2}$, $\bar{b}_{i1}$ are all bounded. Recalling the definition of $\alpha_i$, it can be obtained that $v_i$ is bounded. Furthermore, $x_i$ is bounded as $y_i$ is bounded by Lemma \ref{lemma3}. To this end, we have shown that all the variables contained in the closed-loop system are bounded for $t\in[0,t_f),$ indicating that there is no finite-time escape and $t_f=\infty.$

Step 3: The rest analysis follows the proof of Theorem \ref{th2} to take the time derivatives of $\dot{k}_{i1}$ and $\dot{\hat{\theta}}_{i1}$ to show that there are uniformly continuous. Then, take the integrations of them over $[0,\infty)$ to prove that their integrations are bounded. With the above conclusions in mind, by Barbalat's lemma, $\lim_{t\rightarrow \infty}\dot{k}_{i1}=0$ and $\lim_{t\rightarrow \infty}\dot{\hat{\theta}}_{i1}=0,$ from which it can be obtained that $\lim_{t\rightarrow\infty} x_i(t)-y_i(t)=0$ by following the arguments in the proof of  Theorem \ref{theorem2}.
\end{Proof}
\begin{Remark}
The system dynamics considered in \eqref{thri} is similar to the one in  \cite{YETAC98} and the state regulation part is motivated by \cite{YETAC98}. However, different from \cite{YETAC98} that regulates the state to zero, this paper needs to regulate the state to a time-varying reference trajectory ($y_i(t)$ for $i\in\mathcal{N}$), generated by the optimization module.
%Hence, it is clear that the authors in \cite{YETAC98} addressed an adaptive regulation problem, while this paper needs to achieve distributed Nash equilibrium seeking for the considered system.
\end{Remark}

\section{Discussions on Fully Distributed Implementation of the Proposed Algorithms}\label{dis}
In Section \ref{ne_se_no}, the proposed seeking strategies contain a centralized control gain $\delta$, which depends on the players' objective functions and the communication graph. In general, these centralized information can hardly be obtained. Actually, in \cite{YeATsub}, we proposed fully distributed Nash equilibrium seeking strategies by adaptively adjusting the control gains. In the following, we further prove that the adaptive algorithms in \cite{YeATsub} can also be utilized in the proposed algorithms to achieve fully distributed Nash equilibrium seeking in the considered problem.

By the methods in \cite{YeATsub}, one can replace \eqref{eq3}-\eqref{eq4} in the proposed algorithms with
\begin{equation}\label{eqqq3}
\begin{aligned}
\dot{y}_i=&-\nabla_if_i(\mathbf{z}_i)\\
\dot{z}_{ij}=&-\delta_{ij}\left(\sum_{k=1}^Na_{ik}(z_{ij}-z_{kj})+a_{ij}(z_{ij}-y_j)\right)\\
\dot{\delta}_{ij}=&\left(\sum_{k=1}^Na_{ik}(z_{ij}-z_{kj})+a_{ij}(z_{ij}-y_j)\right)^2,
\end{aligned}
\end{equation}
for $i\in\mathcal{N}.$

Then, the following result can be obtained.
\begin{Lemma}\label{lemma44}
Suppose that Assumptions \ref{ASS1}-\ref{ass_3} are satisfied. Then, with the strategy in \eqref{eqqq3}, the following conclusions can be obtained:
\begin{itemize}
  \item For each $i,j\in\mathcal{N},$ $y_i(t)$, $z_{ij}(t)$ and $\delta_{ij}(t)$ are bounded for $t\in[0,\infty).$
 \item For each $i\in\mathcal{N},$  $\dot{y}_i(t)$ is square integrable over $t\in(0,\infty)$, i.e., $\int_{0}^{\infty} \dot{y}_i^2(s)ds\leq c_i$ for some positive constant $c_i.$
\end{itemize}
\end{Lemma}
 \begin{Proof}
 Following the proof of \cite{YeATsub} to define
 $V=\mathbf{e}^TM\mathbf{e}+\frac{1}{2}(\mathbf{y}-\mathbf{x}^*)^T(\mathbf{y}-\mathbf{x}^*)+\sum_{i=1}^N\sum_{j=1}^N (\theta_{ij}-\theta_{ij}^*)^2,$
where $\theta_{ij}^*>\frac{8m||M||\sqrt{N}\max_{i\in\mathcal{V}}\{l_i\}+(2||M||N+\max_{i\in\mathcal{V}}\{l_i\})^2}{8m\lambda_{min}(MM)},
$ $\mathbf{e}=[z_{11}-y_1,z_{12}-y_2,\cdots,z_{1N}-y_N,z_{21}-y_1,\cdots,z_{NN}-y_N]^T,$ $\mathbf{y}=[y_1,y_2,\cdots,y_N]^T,$  $M=\mathcal{L}\otimes I_{N\times N}+\mathcal{A}_0$ and $\mathcal{A}_0$ is a diagonal matrix with its elements being $a_{ij}.$

Then, it follows from \cite{YeATsub} that
\begin{equation}\label{time_d}
\dot{V}\leq -a||\mathbf{E}||^2,
\end{equation}
where $a>0$ and $\mathbf{E}=[(\mathbf{y}-\mathbf{x}^*)^T,\mathbf{e}^T]^T,$ from which it can be obtained that
 for each $i,j\in\mathcal{N},$ $y_i$, $z_{ij}$ and $\delta_{ij}$ are bounded for $t\in[0,\infty).$

Moreover,
\begin{equation}
\begin{aligned}
\int_{0}^{\infty} \dot{y}_i^2(s)ds\leq &\int_{0}^{\infty} |\nabla_if_i(\mathbf{z}_i(s))-\nabla_if_i(\mathbf{x}^*) |^2ds\\
\leq & l_i^2\int_{0}^{\infty} ||\mathbf{z}_i(s)-\mathbf{x}^*||^2ds.
\end{aligned}
\end{equation}

Taking integration on both sides of \eqref{time_d}, we obtain that
\begin{equation}\label{time_d1}
\int_0^{\infty}\dot{V}(s)ds\leq -a\int_0^{\infty}||\mathbf{E}(s)||^2ds,
\end{equation}
by which
\begin{equation}
V(\infty)+a\int_0^{\infty}||\mathbf{E}||^2ds\leq V(0).
\end{equation}
By further noticing that
\begin{equation}
\int_{0}^{\infty} ||\mathbf{z}_i(s)-\mathbf{x}^*||^2ds \leq \int_0^{\infty}||\mathbf{E}||^2ds,
\end{equation}
we obtain that
\begin{equation}
V(\infty)+a\int_{0}^{\infty} ||\mathbf{z}_i(s)-\mathbf{x}^*||^2ds \leq V(0).
\end{equation}
Hence
\begin{equation}
\int_{0}^{\infty} \dot{y}_i^2(s)ds\leq \frac{(V(0)-V(\infty))l_i^2}{a},
\end{equation}
thus arriving at the second conclusion.
 \end{Proof}

With the results in Lemma \ref{lemma44}, we can achieve the fully distributed implementations of the proposed algorithms, which is stated in the following theorem.
\begin{Theorem}\label{thee}
Suppose that Assumptions \ref{ASS1}-\ref{Assum_4} are satisfied. Then, for the system considered in \eqref{dyna_1} with the control input in \eqref{eq1}-\eqref{eq2}, where $y_i$ is generated by \eqref{eqqq3}. Then,
\begin{equation}
\lim_{t\rightarrow\infty}||\mathbf{x}(t)-\mathbf{x}^*||=0,
\end{equation}
 and all the other variables stay bounded.
\end{Theorem}

It's worth mentioning that for systems considered in \eqref{dyna_2}/\eqref{dyna_3} and the proposed control inputs designed for the corresponding systems, one can replace $y_i$ therein with the one generated by \eqref{eqqq3} to achieve fully distributed implementations of the proposed algorithms. Note that we only present the results for the system \eqref{dyna_1} and omit the rest to avoid any repetitions in this paper.
\begin{Remark}
In this section, we only provide an example to illustrate the fully distributed implementations of the proposed algorithms. However, it is worth noting that the proposed algorithms actually provide a general framework to deal with games in systems with unknown control directions. That is, one may utilize other alternative approaches that result in square integrable $\dot{y}_i$ and bounded state variables as well as their time derivatives to achieve fully distributed Nash equilibrium seeking for systems with unknown controls.
\end{Remark}
\begin{Remark}
The modular design in this paper is motivated by \cite{YETACarXiv}\cite{YEAT16}. Though it was required that the controls should be bounded in \cite{YETACarXiv}, the control directions were supposed to be known. Moreover, \cite{YEAT16} designed an extremum seeker through robust state regulation and numerical optimization, in which the control directions are also considered to be known. Different from \cite{Tang} that considered distributed optimization problems with unknown control directions, this paper addresses Nash equilibrium seeking problems with both unknown control directions and parametric uncertainties. In particular, the existence of multiple unknown control directions and uncertain parameters is addressed. Though we only investigate first-order and second-order systems analytically in this paper, we believe that under the proposed framework, it is not challenging to extend the current results to high-order systems by backstepping techniques.
\end{Remark}
Though for presentation simplicity, we suppose that $x_i\in\mathbb{R},$ it should be noted that the presented results can be directly adapted to deal with games in which the players are of multiple heterogeneous dimensions.  In the subsequent section, an example in which $x_i\in\mathbb{R}^2$ for $i\in\mathcal{N}$ will be numerically studied.

\section{A Numerical Example}\label{P_1_exa}
\begin{figure}
  \centering
  % Requires \usepackage{graphicx}
  \includegraphics[scale=0.3]{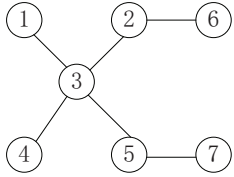}\\
  \caption{The communication graph among the players.}\label{ff1}
\end{figure}

In this section, we consider the connectivity control game among a network of $7$ mobile sensors considered in \cite{YEcyber2020}. The objective function of player $i$ engaged in the game is defined as
\begin{equation}
F_i(\mathbf{x})=h_i(x_i)+l_i(\mathbf{x}),
\end{equation}
where $x_i=[x_{i1},x_{i2}]^T\in \mathbb{R}^2$ and
\begin{equation}
h_i(x_i)=x^T_im_{ii}x_i+x^T_im_i+i^2,
\end{equation}
in which $m_{ii}=diag\{2i,i\},m_i=[i,2i]^T.
$ Moreover,
$l_1(\mathbf{x})={\lVert x_1-x_2\lVert}^2$, $l_2(\mathbf{x})={\lVert x_2-x_3\lVert}^2,$ $l_3(\mathbf{x})={\lVert x_3-x_1\lVert}^2,$
$l_4(\mathbf{x})={\lVert x_4-x_3\lVert}^2$,$l_5(\mathbf{x})={\lVert x_5-x_1\lVert}^2+{\lVert x_5-x_6\lVert}^2$,$l_6(\mathbf{x})={\lVert x_6-x_3\lVert}^2+{\lVert x_6-x_1\lVert}^2$ and $l_7(\mathbf{x})={\lVert x_7-x_2\lVert}^2$. It can be calculated that the Nash equilibrium of the game is  $x^*_{i1}=-\frac{1}{4}$ and $x^*_{i2}=-1$ for $i\in \{1,2,\cdots,7\}$. In the simulation, the undirected and connected communication graph is plotted in
Fig. \ref{ff1}. In the following, games with dynamics in \eqref{dyna_1}, \eqref{dyna_2} and \eqref{dyna_3} will be numerically explored, successively.

\subsection{Distributed Nash equilibrium seeking for first-order systems}\label{fireq}
\begin{figure}
  \centering
  % Requires \usepackage{graphicx}
  \includegraphics[scale=0.5]{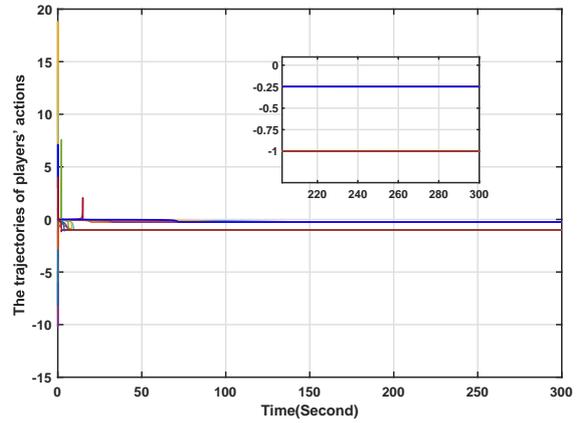}\\
  \caption{The trajectories of $x_i(t)$ for $i\in\{1,2,\cdots,7\}$ generated by \eqref{eq1}-\eqref{eq4}.}\label{ffig3}
\end{figure}

\begin{figure}
  \centering
  % Requires \usepackage{graphicx}
  \includegraphics[scale=0.5]{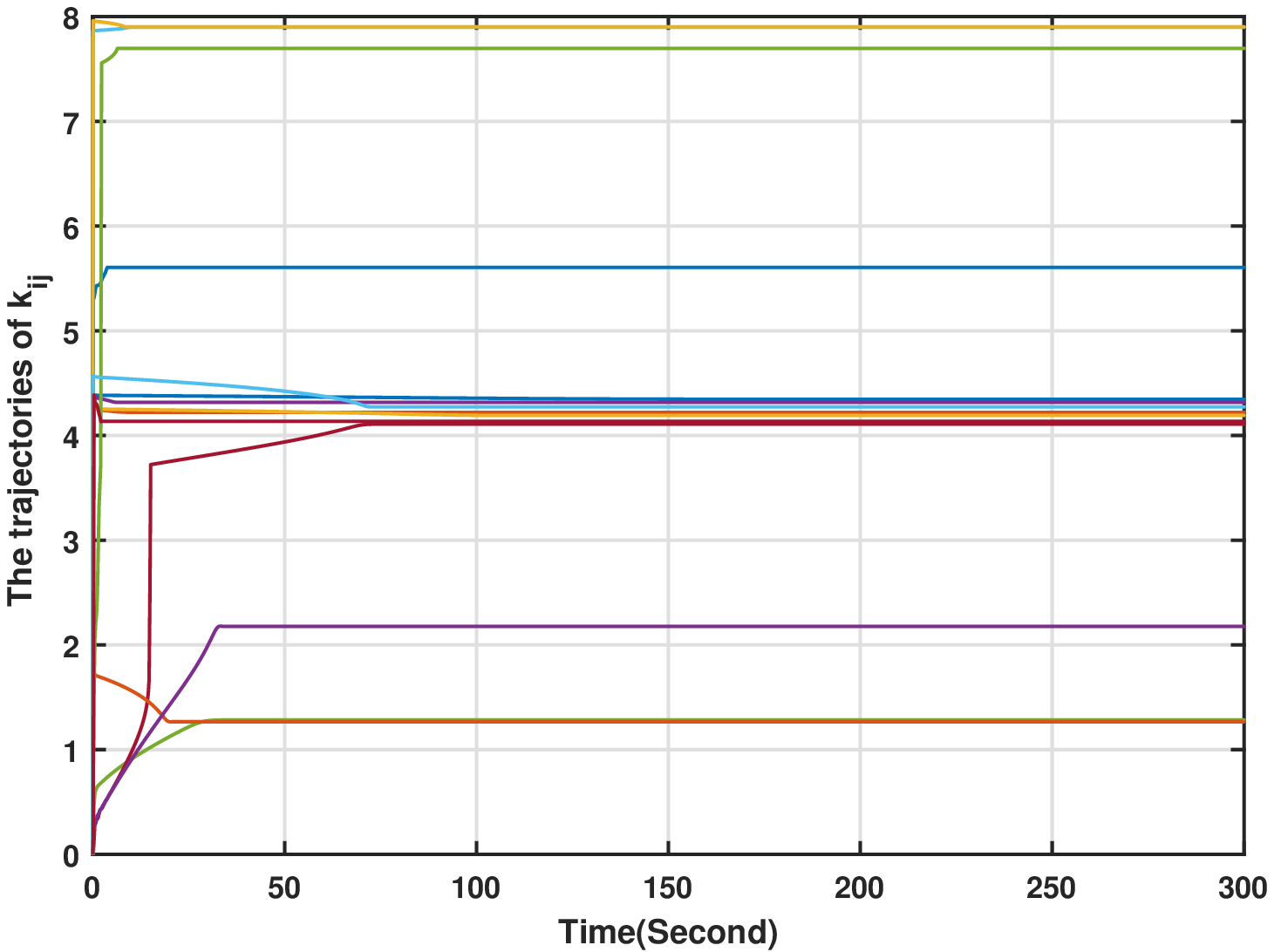}\\
  \caption{The trajectories of $k_{ij}$  for $i\in \{1,2,\cdots,7\}, j\in \{1,2\}$ generated by \eqref{eq1}-\eqref{eq4}.}\label{ffig4}
\end{figure}

\begin{figure}
  \centering
  % Requires \usepackage{graphicx}
  \includegraphics[scale=0.52]{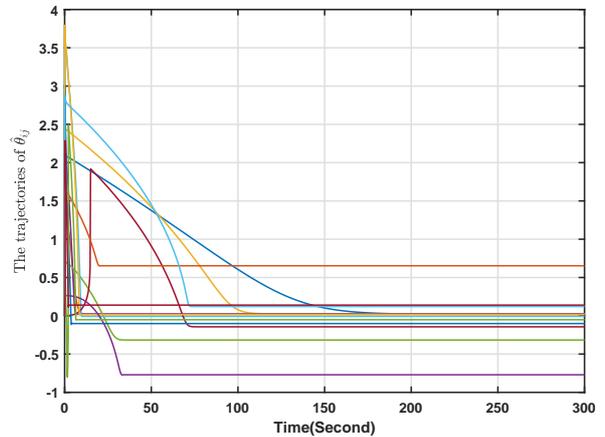}\\
  \caption{The trajectories of $\hat{\theta}_{ij}$ generated by \eqref{eq1}-\eqref{eq4} for $i\in \{1,2,\cdots,7\}, j\in \{1,2\}$.}\label{ffig5}
\end{figure}

In this section, we simulate first-order systems in \eqref{dyna_1}, where the control input is designed in \eqref{eq1}-\eqref{eq4}. Note that as $\mathbf{x}_i\in\mathbb{R}^2,$ $b_i\in\mathbb{R}^2\times \mathbb{R}^2.$ In the simulation, $b_1=diag\{3,3\},
b_2=diag\{5,5\},$  $b_3=diag\{-2,-2\},$ $b_4=diag\{1,2\},$ $b_5=diag\{-3,-3\},$ $b_6=diag\{-1,-1\}$ and $b_7=diag\{2,2\}.$
Moreover, $\phi_i=ix_{i}.$ Let  $\mathbf{x}(0)=[-5,3,-4,-6,1,8,0,-8,-1,10,1,2,3,0]^T,$ and the initial values for all the other variables in \eqref{eq1}-\eqref{eq4} be zero. Then, generated by \eqref{eq1}-\eqref{eq4}, the players' action trajectories $x_i(t)$ for $i\in\{1,2,\cdots,7\}$ are plotted in Fig. \ref{ffig3}, from which it is clear that the players' action trajectories converge to the Nash equilibrium asymptotically. Moreover, Figs. \ref{ffig4}-\ref{ffig5} illustrate the trajectories of $k_{ij}(t)$ and $\hat{\theta}_{ij}(t)$ for all $i\in \{1,2,\cdots,7\}, j\in \{1,2\}$, respectively. From Figs. \ref{ffig4}-\ref{ffig5}, it can be seen that these variables stay bounded. Therefore, Theorem \ref{theorem2} is numerically validated.

\subsection{Distributed Nash equilibrium seeking for second-order systems}\label{sec}

In this section, we simulate the system in \eqref{dyna_2}, where the control input is designed in  \eqref{eqq1}-\eqref{eqq4}. In the simulation, $b_i$, $\phi_i(x_i)$ and $\mathbf{x}(0)$ follow those in Section \ref{fireq} and the initial values for all the other variables are zero.

\begin{figure}
  \centering
  % Requires \usepackage{graphicx}
  \includegraphics[scale=0.55]{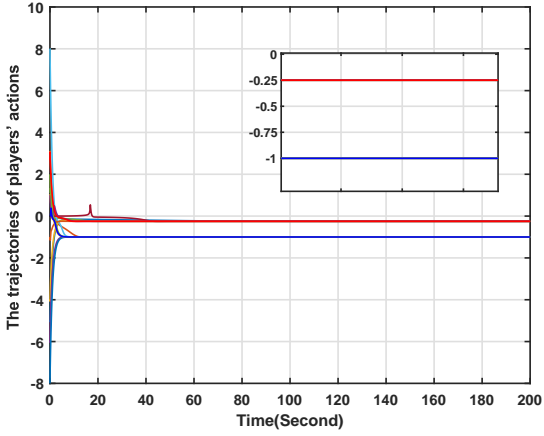}\\
  \caption{The trajectories of $x_i(t)$ for $i\in\{1,2,\cdots,7\}$ generated by \eqref{eqq1}-\eqref{eqq4}.}\label{figg5}
\end{figure}

\begin{figure}
  \centering
  % Requires \usepackage{graphicx}
  \includegraphics[scale=0.5]{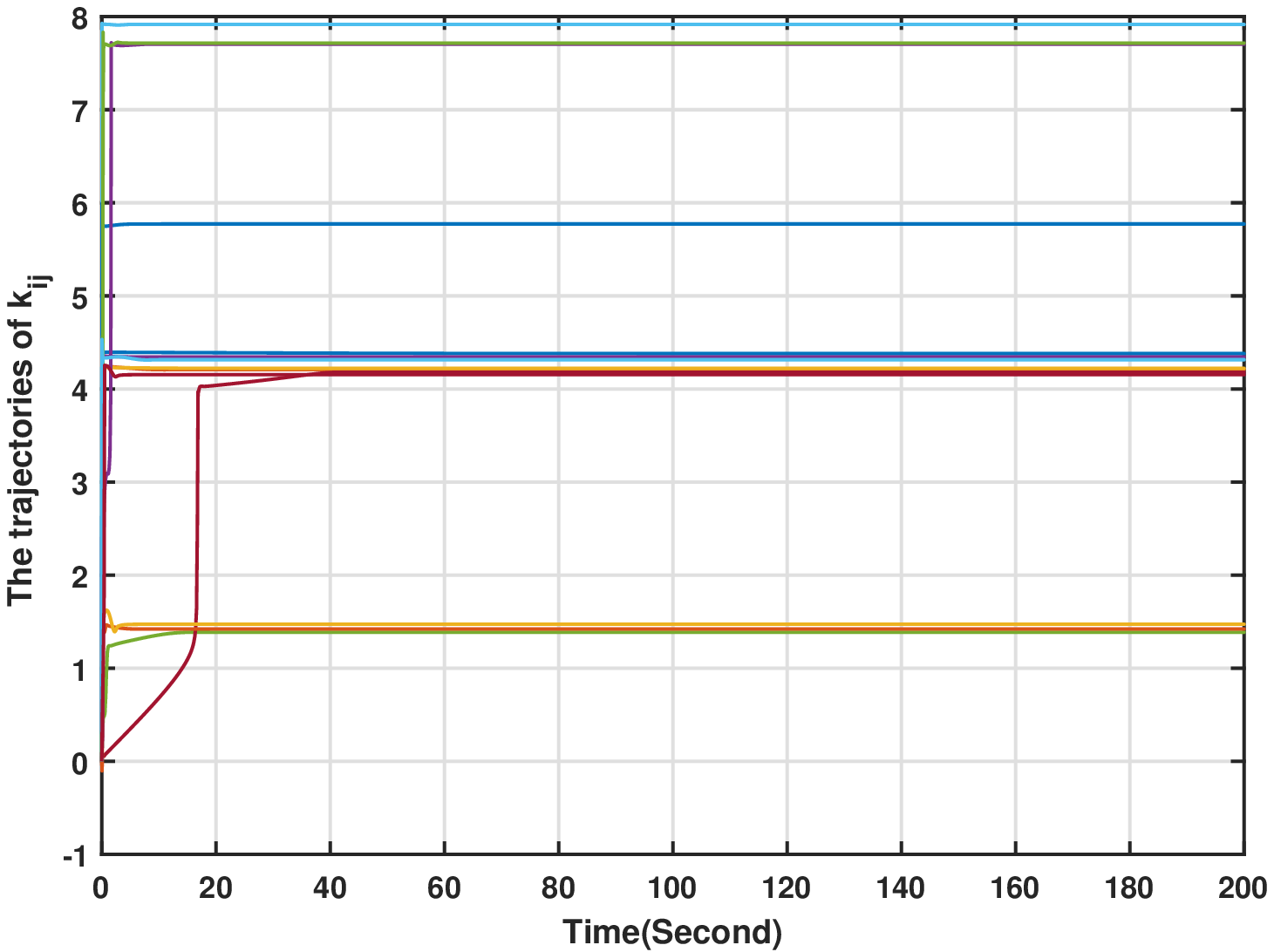}\\
  \caption{The trajectories of $k_{ij}$ for $i\in \{1,2,\cdots,7\}, j\in \{1,2\}$  generated by \eqref{eqq1}-\eqref{eqq4}.}\label{figg6}
\end{figure}

\begin{figure}
  \centering
  % Requires \usepackage{graphicx}
  \includegraphics[scale=0.5]{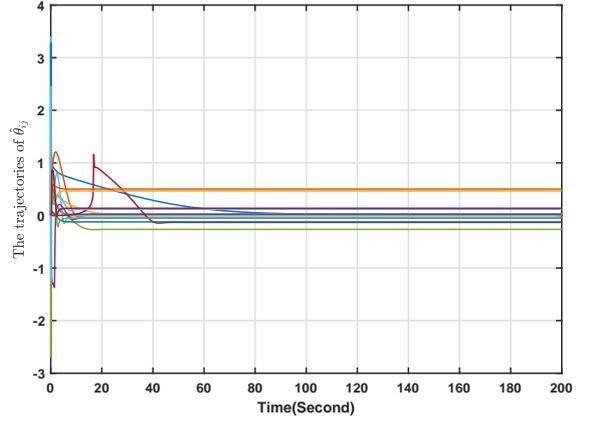}\\
  \caption{The trajectories of $\hat{\theta}_{ij}$ for $i\in \{1,2,\cdots,7\}, j\in \{1,2\}$  generated by  generated by \eqref{eqq1}-\eqref{eqq4}.}\label{figg7}
\end{figure}

\begin{figure}
  \centering
  % Requires \usepackage{graphicx}
  \includegraphics[scale=0.5]{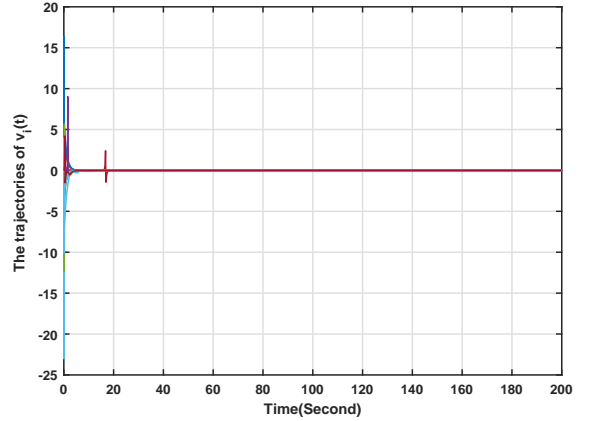}\\
  \caption{The trajectories of $v_i(t)$ for $i\in\{1,2,\cdots,7\}$ generated by \eqref{eqq1}-\eqref{eqq4}.}\label{figg8}
\end{figure}

The players' action trajectories $x_i(t)$ for $i\in\{1,2,\cdots,7\}$ generated by \eqref{eqq1}-\eqref{eqq4} are depicted in Fig. \ref{figg5}, from which it can be seen that the players' actions converge the actual Nash equilibrium of the game. In addition, $k_{ij}(t)$ and $\hat{\theta}_{ij}(t)$ for all $i\in \{1,2,\cdots,7\}, j\in \{1,2\}$ are given in Figs. \ref{figg6}-\ref{figg7}. From Figs. \ref{figg6}-\ref{figg7}, we can conclude that $k_{ij}(t)$ and $\hat{\theta}_{ij}(t)$ for all $i\in \{1,2,\cdots,7\}, j\in \{1,2\}$ stay bounded. Furthermore, Fig. \ref{figg8} demonstrates that $v_i(t)$ for all $i\in\{1,2,\cdots,7\}$ decay to zero, which is aligned with the results in Theorem \ref{th2}. To this end, the conclusions in Theorem \ref{th2} have been numerically verified.


\subsection{Distributed Nash equilibrium seeking for more general second-order uncertain nonlinear systems}
\begin{figure}
  \centering
  % Requires \usepackage{graphicx}
  \includegraphics[scale=0.5]{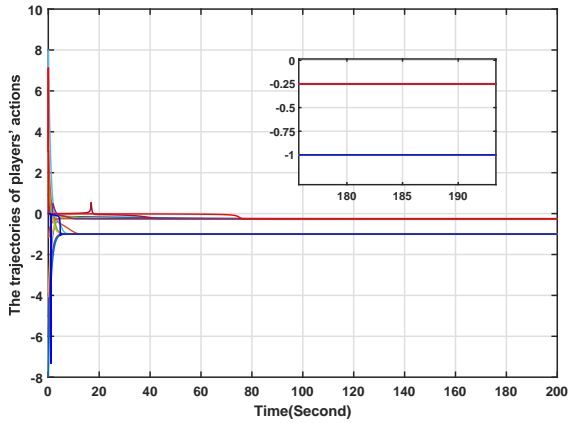}\\
  \caption{The trajectories of $x_i(t)$ for $i\in\{1,2,\cdots,7\}$ with player $7$'s control strategy being \eqref{eqq13}-\eqref{eqq15} and the rest of the players' control strategy being \eqref{eqq1}-\eqref{eqq4}.}\label{Th3x}
\end{figure}

\begin{figure}
  \centering
  % Requires \usepackage{graphicx}
  \includegraphics[scale=0.5]{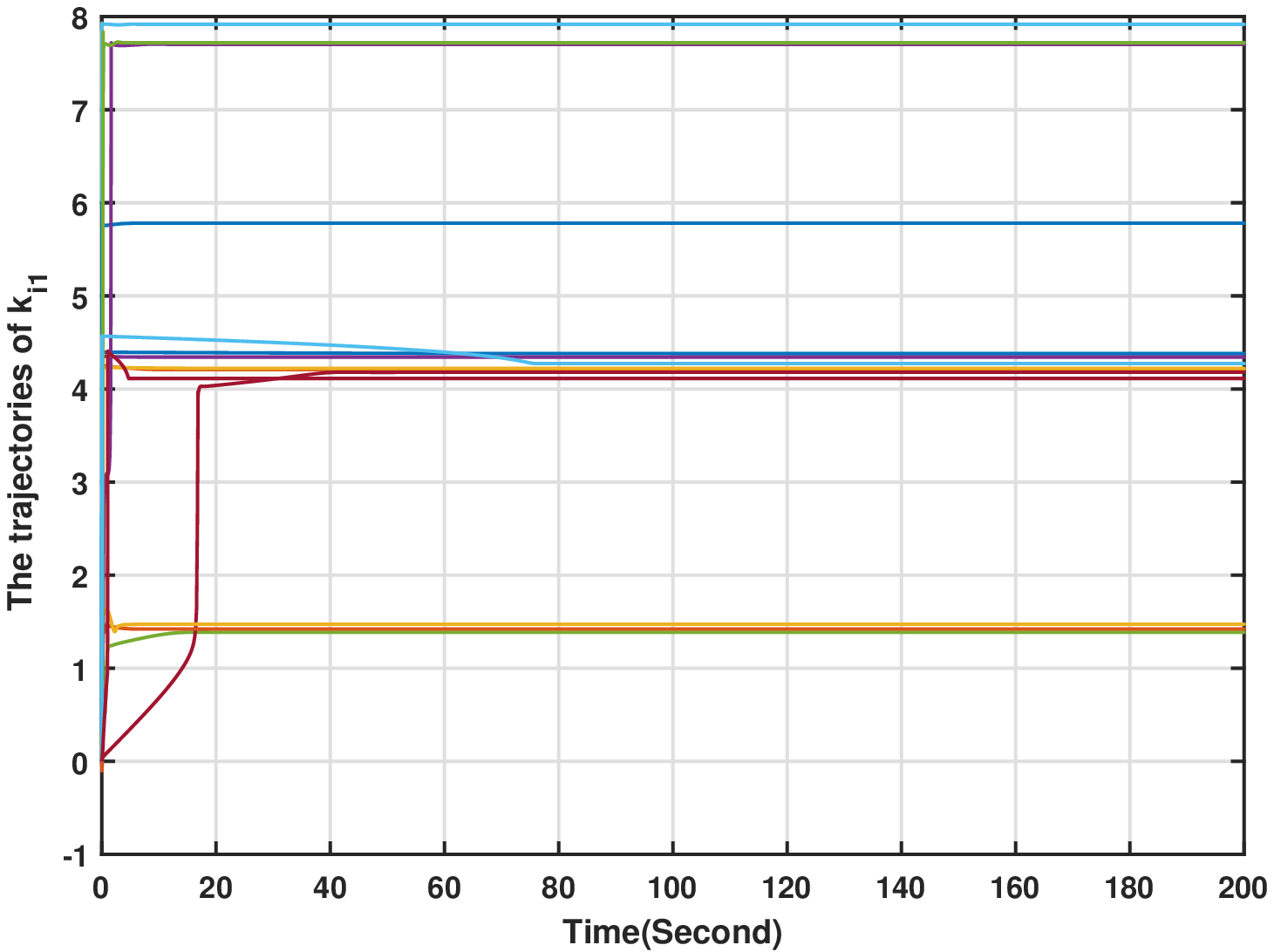}\\
  \caption{The trajectories of $k_{i1}(t)$ for $i\in \{1,2,\cdots,7\}$ with player $7$'s control  strategy being \eqref{eqq13}-\eqref{eqq15} and the rest of the players' control strategy being \eqref{eqq1}-\eqref{eqq4}.}\label{Th3ki1}
\end{figure}

\begin{figure}
  \centering
  % Requires \usepackage{graphicx}
  \includegraphics[scale=0.52]{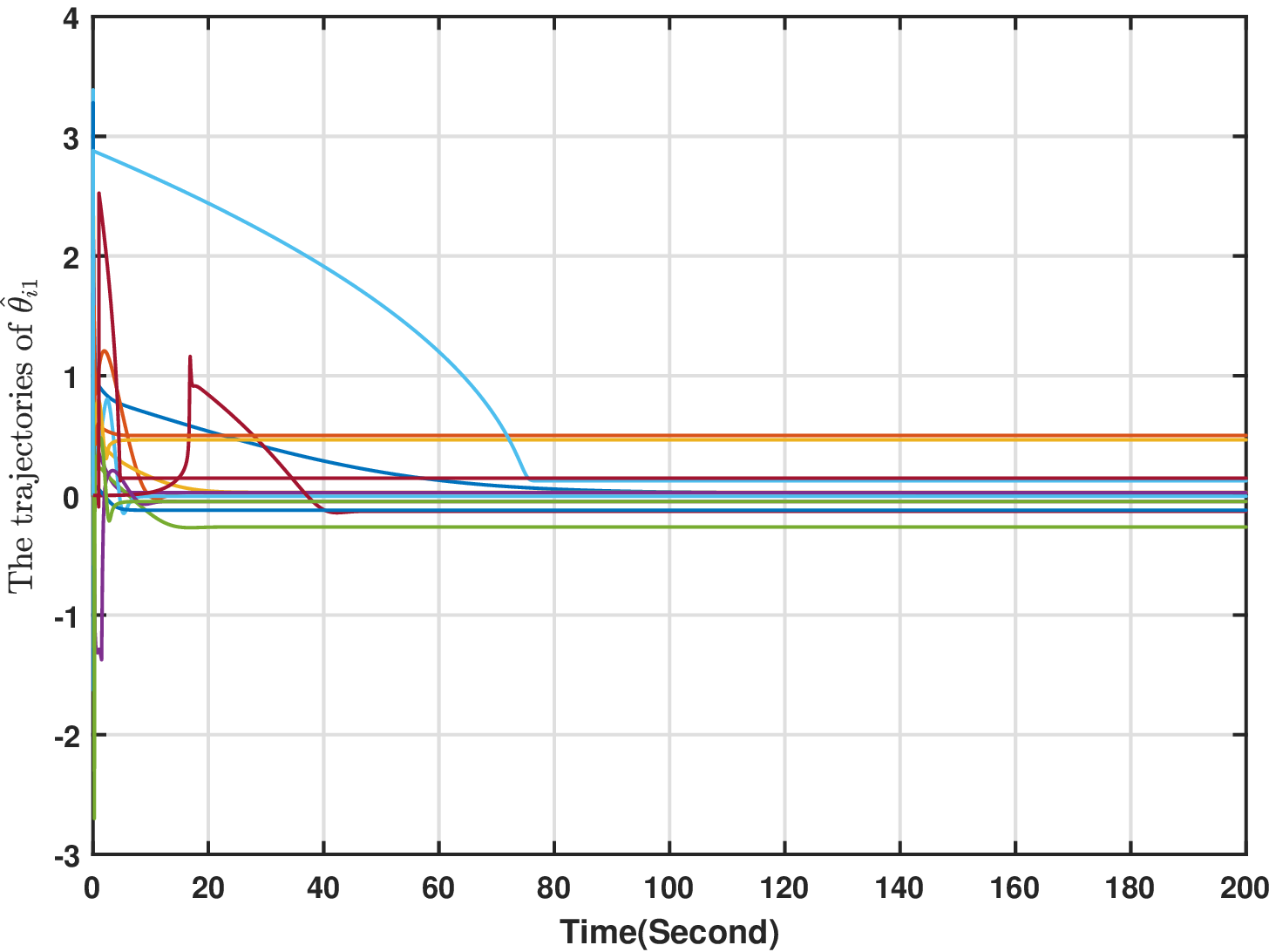}\\
  \caption{The trajectories of $\hat{\theta}_{i1}(t)$  for $i\in \{1,2,\cdots,7\}$ with player $7$'s control  strategy being \eqref{eqq13}-\eqref{eqq15} and the rest of the players' control strategy being \eqref{eqq1}-\eqref{eqq4}.}\label{Th3theati1}
\end{figure}

\begin{figure}
  \centering
  % Requires \usepackage{graphicx}
  \includegraphics[scale=0.5]{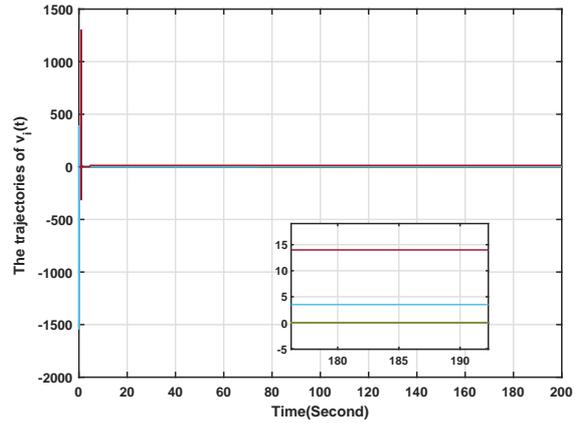}\\
  \caption{The trajectories of $v_i(t)$ for $i\in\{1,2,\cdots,7\}$ with player $7$'s control input being \eqref{eqq13}-\eqref{eqq15} and the rest of the players' control input being \eqref{eqq1}-\eqref{eqq4}.}\label{Th3v}
\end{figure}

In this section, we numerically verify the control input designed for uncertain nonlinear systems in \eqref{dyna_3}. To illustrate the case, we suppose that the action of player $7$ is governed by \eqref{dyna_3}, and its control input is given by  \eqref{eqq13}-\eqref{eqq15}, while all the other players' actions are governed by \eqref{dyna_2} with their control inputs being \eqref{eqq1}-\eqref{eqq4}.
For players $1$-$6$, $b_{i}$ and $\phi_i(x_i)$ are chosen to be the same as those in Section \ref{sec}. For player $7$, $b_{71}=b_{72}=diag\{2,2\},$ $\phi_{71}(x_i)=7x_{7}$ and $\phi_{72}(x_i,v_i)=[7x_{72},7v_{72}]^T.$
In addition,  $\mathbf{x}(0)=[-5,3,-4,-6,1,8,0,-8,-1,10,1,2,3,0]^T$ and the initial conditions for all the other variables are zero. Generated by the proposed methods, the players' actions $x_i(t)$ for $i\in\{1,2,\cdots,7\}$ are shown in Fig. \ref{Th3x}, from which we see that the players' actions converge to the Nash equilibrium. Moreover,
$k_{i1}(t)$ and $\hat{\theta}_{i1}(t)$ for $i\in \{1,2,\cdots,7\}$ are plotted in Figs. \ref{Th3ki1}-\ref{Th3theati1}. Figs. \ref{Th3ki1}-\ref{Th3theati1} show that
$k_{i1}(t)$ and $\hat{\theta}_{i1}(t)$ stay bounded.  Moreover, the evolution of $v_i(t)$ is shown in Fig. \ref{Th3v}, which shows that $v_i(t)$ for all $i\in\{1,2,\cdots,7\}$ are also bounded. Hence, the effectiveness of the method in \eqref{eqq13}-\eqref{eqq15} is also verified.

\subsection{Fully distributed implementations of the proposed algorithms}
\begin{figure}
  \centering
  % Requires \usepackage{graphicx}
  \includegraphics[scale=0.5]{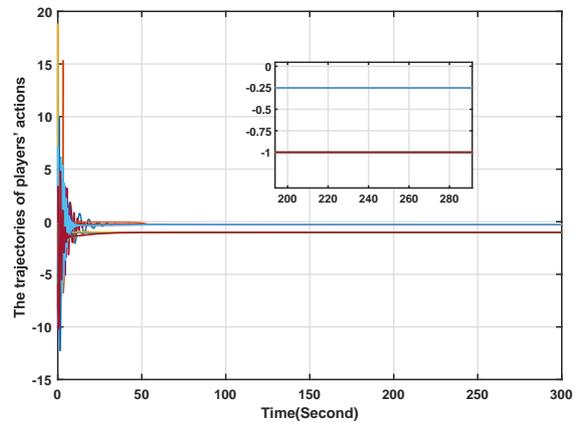}\\
  \caption{The trajectories of $x_i(t)$ for $i\in\{1,2,\cdots,7\}$ for the system considered in \eqref{dyna_1} with the control input in \eqref{eq1}-\eqref{eq2}, where $y_i$ is generated by \eqref{eqqq3}.}\label{xf}
\end{figure}

\begin{figure}
  \centering
  % Requires \usepackage{graphicx}
  \includegraphics[scale=0.5]{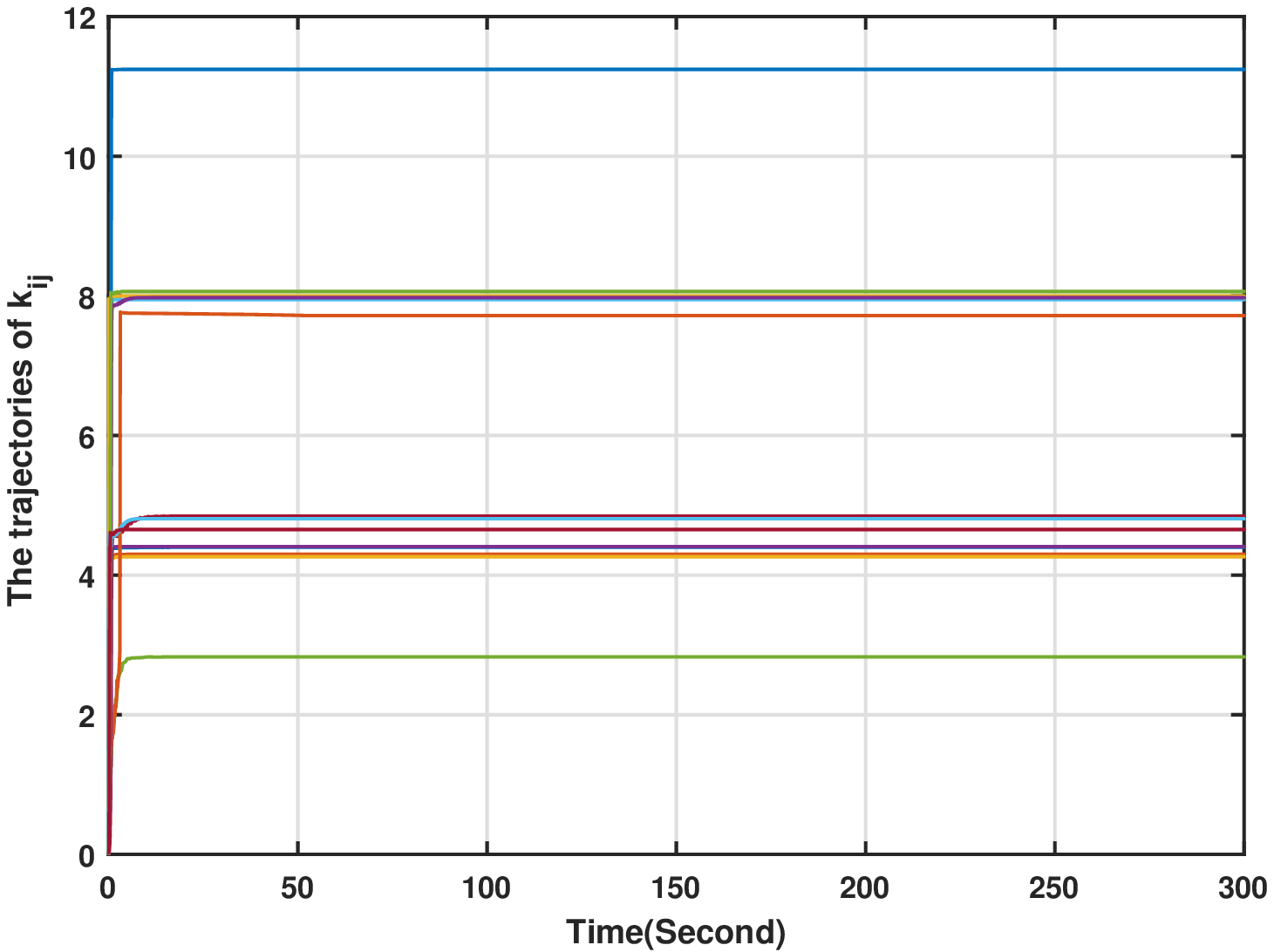}\\
  \caption{The trajectories of $k_{i}(t)$ for $i\in \{1,2,\cdots,7\}$ for the system considered in \eqref{dyna_1} with the control input in \eqref{eq1}-\eqref{eq2}, where $y_i$ is generated by \eqref{eqqq3}.}\label{kf}
\end{figure}

\begin{figure}
  \centering
  % Requires \usepackage{graphicx}
  \includegraphics[scale=0.52]{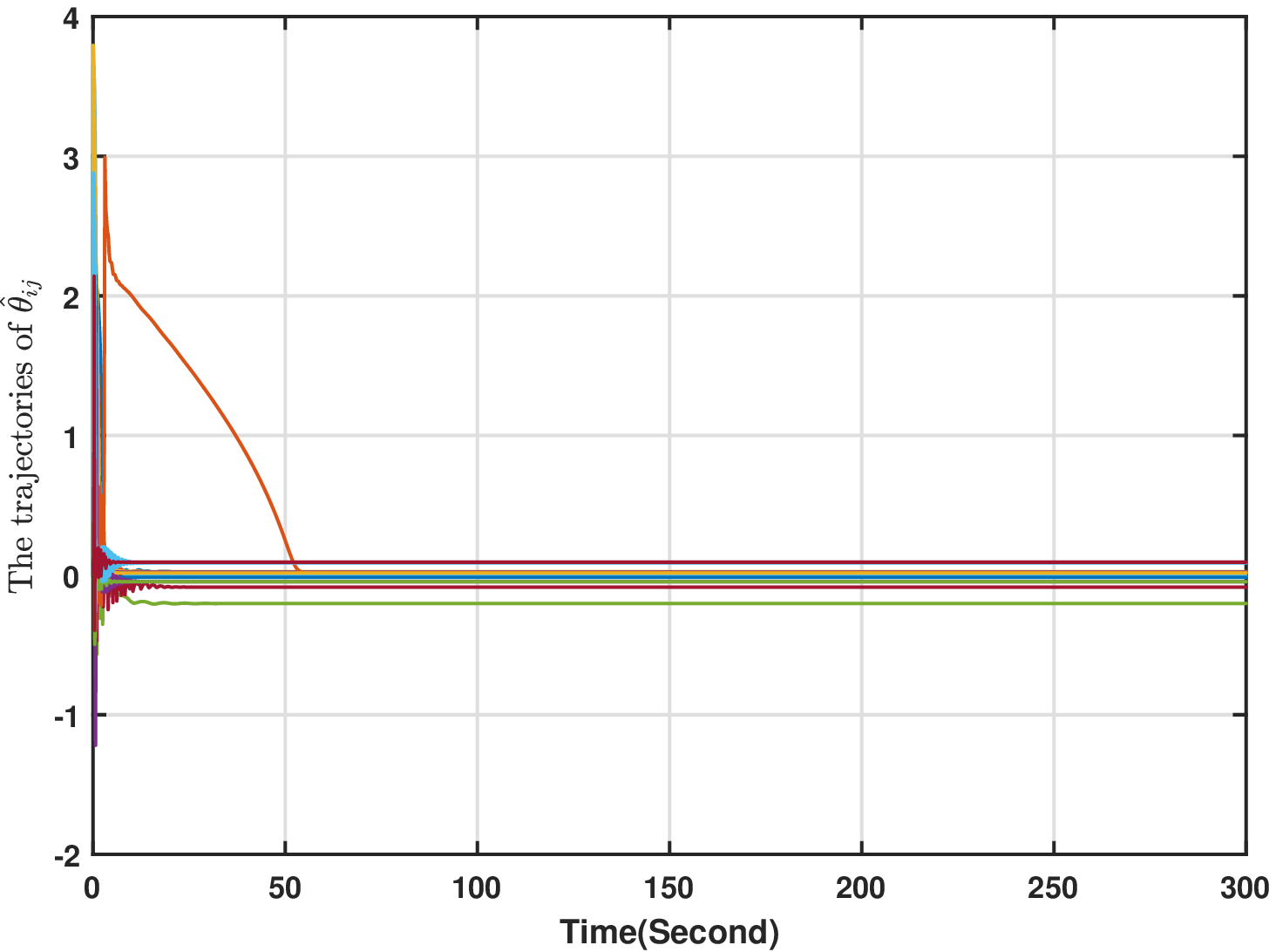}\\
  \caption{The trajectories of $\hat{\theta}_{ij}(t)$  for $i\in \{1,2,\cdots,7\}$ for the system considered in \eqref{dyna_1} with the control input in \eqref{eq1}-\eqref{eq2}, where $y_i$ is generated by \eqref{eqqq3}.}\label{theatf}
\end{figure}

\begin{figure}
  \centering
  % Requires \usepackage{graphicx}
  \includegraphics[scale=0.52]{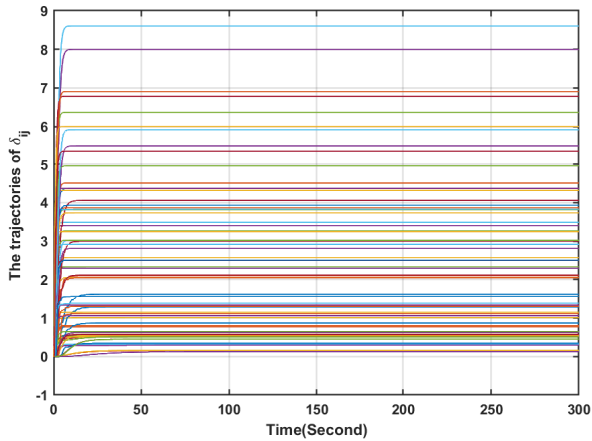}\\
  \caption{The trajectories of $\delta_{ij}(t)$ for $i\in\{1,2,\cdots,7\},j\in\{1,2\}$  for the system considered in \eqref{dyna_1} with the control input in \eqref{eq1}-\eqref{eq2}, where $y_i$ is generated by \eqref{eqqq3}.}\label{delta_f}
\end{figure}

To verify the fully distributed implementations of the proposed methods, we take first-order systems as an example. The simulation setting of this section follows that of Section \ref{fireq} and $\delta_{ij}(0)=0$ for all $i\in\{1,2,\cdots,7\},j\in\{1,2\}$. The simulation results  for the system considered in \eqref{dyna_1} with the control input in \eqref{eq1}-\eqref{eq2}, where $y_i$ is generated by \eqref{eqqq3} are given in Figs. \ref{xf}-\ref{delta_f}. Fig. \ref{xf} plots the players' actions from which we see that the players' actions can converge to the Nash equilibrium. Moreover, Figs. \ref{kf}-\ref{delta_f} plot $k_{i}(t),$ $\hat{\theta}_{i}(t)$ and $\delta_{ij}(t)$, respectively, from which it is clear that they stay bounded. Therefore, the results in Theorem \ref{thee} is numerically verified. Note that compared with the simulation in Section 6.1, there is no centralized control gain in \eqref{eqqq3}, thus verifying the effectiveness of the distributively implemented algorithms.


\section{Conclusion}\label{P_1_CONC}
This paper considers distributed Nash equilibrium seeking for games in which the players' actions are subject to both unknown control directions and parametric uncertainties. First-order systems and second-order systems are addressed successively. To cope with the un-availability of control directions, a Nussbaum function is adopted. Moreover, the parametric uncertainties are addressed by adaptive laws. Together with an optimization module, a state regulation module is included in the seeking strategy.  Based on the Barbalat's lemma, it is proven that the players' actions can be driven to the Nash equilibrium. Lastly, the fully distributed implementations of the proposed algorithms are discussed. It is shown that the adaptive techniques can be employed to achieve the equilibrium seeking in a fully distributed way.

%\ \ \\
%\begin{wrapfigure}{l}{25mm}
%    \includegraphics[width=1in,height=1.25in,clip,keepaspectratio]{mjye1.eps}
%  \end{wrapfigure}\par
%  \textbf{Maojiao Ye} received the B.Eng. degree in Automation from the University of Electronic Science and Technology of China, Sichuan, China, in 2012 and the Ph.D. degree from Nanyang Technological University, Singapore, in 2016. She is currently a Professor with the School of Automaton, Nanjing University of Science and Technology. Prior to her current position, she was a research fellow in the School of Electrical and Electronic Engineering at Nanyang Technological University from 2016-2017.
%
%  Dr. Ye was a recipient of Guan Zhao-Zhi Award in the 36th Chinese Control Conference 2017 (first author) and a recipient of the Best Paper Award in the 15th IEEE International Conference on Control and Automation 2019 (sole author). Her research interests include game theory, distributed optimization, extremum seeking control and demand response in smart grids.


\begin{thebibliography}{99}
\bibitem{YEcyber2020}M. Ye, ``Distributed robust seeking of Nash equilibrium for networked games: an extended state observer-based approach," \emph{IEEE Transactions on Cybernetics}, accepted, published online, DOI: 10.1109/TCYB.2020.2989755.
\bibitem{Yecyber17}M. Ye and G. Hu, ``Game design and analysis for price-based demand response: an aggregate game approach," \emph{IEEE Transactions on Cybernetics,} vol. 47, no. 3, pp. 720-730, 2017.
\bibitem{YETAC19}M. Ye, G. Hu, F. Lewis, L. Xie, ``A unified strategy for solution seeking in graphical $N$-coalition noncooperative games," \emph{IEEE Transactions on Automatic Control}, vol. 64, no. 11, pp. 4645-4652, 2019.
\bibitem{MAIJRNC19}J. Ma, M. Ye, Y. Zheng, Y. Zhu, ``Consensus analysis of hybrid multiagent systems: a game-theoretic approach," \emph{International Journal of Robust and Nonlinear Control,} vol.29, pp. 1840-1853, 2019.
\bibitem{YETAC17} M. Ye and G. Hu, ``Distributed Nash equilibrium seeking by a consensus based approach," \emph{IEEE Transactions on Automatic Control,} pp. 4811-4818, vol. 62, no. 9, 2017.
\bibitem{JiangRobo02}P. Jiang, P. Woo, R. Unbehauen, ``Iterative learning control for manipulator trajectory tracking without any control singularity," \emph{Robotica,} vol. 20, no. 2, pp. 149-158, 2002.
\bibitem{DuJOE}J. Du, C. Guo, S. Yu, Y. Zhao, ``Adaptive autopilot design of time-varying uncertain ships with completely unknown control coefficients," \emph{IEEE Journal of Oceanic Engineering,} vol. 32, no. 2, pp. 346-352, 2007.
\bibitem{Nussbaum83}R. Nussbaum, ``Some remarks on a conjecture in parameter adaptive control," \emph{System and Control Letters,} vol. 3, no. 5, pp. 243-246, 1983.
\bibitem{BunJFI16}X. Bun, D. Wei, X. Wu, J. Huang, ``Guaranteeing preselected tracking quality for air-breathing hypersonic non-affine models with an unknown control direction via concise neural control," \emph{Journal of the Franklin Institute,} vol. 353, no. 13, pp. 3207-3232, 2016.
\bibitem{Wang2020}L. Wang, W. Deng, J. Liu and R. Mei,  ``Adaptive sliding mode trajectory tracking control of quadrotor UAV with unknown control direction,"  In: R. Wang, Z. Chen, W. Zhang, Q. Zhu(eds) Proceedings of the 11th International Conference on Modelling, Identification and Control, Lecture Notes in Electrical Engineering, vol. 582. Springer, Singapore.
\bibitem{PengSCL14}J. Peng, X. Ye, ``Cooperative control of multiple heterogeneous agents with unknown high-frequency-gain signs," \emph{Systems and Control Letters,} vol. 68, pp. 51-56, 2014.
\bibitem{YETAC98}X. Ye and J. Jiang, ``Adaptive nonlinear design without a priori knowledge of control directions," \emph{IEEE Transactions on Automatic Control}, vol. 43, no. 11, pp. 1617-1621, 1998.
\bibitem{ChenAT19}Z. Chen, ``Nussbaum functions in adaptive control with time-varying unknown control coefficients," \emph{Automatica}, vol. 102, pp. 72-79, 2019.
\bibitem{Khailil02} H. Khailil, \emph{Nonlinear Systems,} Upper Saddle River, NJ: Prentice Hall, 2002.
\bibitem{ScheinkerTAC13}A. Scheinker, M. Krstic, ``Minimum-seeking for CLFs: universal semiglobally stabilizing feedback under unknown control directions," \emph{IEEE Transactions on Automatic Control}, vol. 58, no. 5, pp. 1107-1122, 2013.
\bibitem{ChenTAC17}C. Chen, C. Wen, Z. Liu, K. Xie, Y. Zhang, C. Chen, ``Adaptive consensus of nonlinear multi-agent systems with non-identical partially unknown control directions and bounded modelling errors," \emph{IEEE Transactions on Automatic Control,} vol. 62, no. 9, pp. 4654-4659, 2017.
\bibitem{YangAT08}C. Yang, S. Ge, T. Lee, ``Output feedback adaptive control of a class of nonlinear discrete-time systems with unknown control directions," \emph{Automatica}, vol. 45, no. 1, pp. 270-276, 2008.
\bibitem{GuoTAC17} M. Guo, D. Xu, L. Liu, ``Cooperative output regulation of heterogeneous nonlinear multi-agent systems with unknown control directions," \emph{IEEE Transactions on Automatic Control}, vol. 62, no. 6, pp 3039-3045, 2017.
\bibitem{Tang}Y. Tang, ``Multi-agent optimal consensus with unknown control directions," 	arXiv:2005.10492.
\bibitem{KoshalOR16}J. Koshal, A. Nedic and U. Shanbhag, ``Distributed algorithms
for aggregative games on graphs," \emph{Operations Research,} vol. 64,
pp. 680-704, 2016.
\bibitem{SalehisadaghianiAT16}F. Salehisadaghiani and L. Pavel, ``Distributed Nash equilibrium seeking:
A gossip-based algorithm," \emph{Automatica,} vol. 72, pp. 209-216,
2016.
\bibitem{YEAT2020}M. Ye, G. Hu and S. Xu, ``An extremum seeking-based approach
for Nash equilibrium seeking in $N$-cluster noncooperative games,"
\emph{Automatica,} vol. 114, 108815, 2020.

\bibitem{YEAT18}M. Ye, G. Hu, and F. L. Lewis, ``Nash equilibrium seeking for $N$-coalition
non-cooperative games," \emph{Automatica,} vol. 95, pp. 266-272,
2018.

\bibitem{Ibrahim18}A. Ibrahim, T. Hayakawa, ``Nash equilibrium seeking with secondorder
dynamic agents," \emph{IEEE Conference on Decision and Control,}
pp. 2514-2518, 2018.
\bibitem{YETACarXiv}M. Ye, ``Distributed Nash equilibrium seeking for games in systems
with bounded control inputs," submitted to \emph{IEEE Transactions on Automatic Control,} revised, available online: arXiv:1901.09333, 2019.
\bibitem{Slotine}J. Slotine, W. Li, \emph{Applied Nonlinear Control,} Prentice Hall, Englewood Cliffs, 1991.
\bibitem{MYEcyber18} M. Ye, G. Hu,``Distributed Nash equilibrium seeking in multi-agent games under switching communication topologies," \emph{IEEE Transactions on Cybernetics,} vol. 48, no. 11, pp. 3208-3217, 2018.
\bibitem{YeATsub}M. Ye, G. Hu, ``Adaptive approaches for fully distributed Nash equilibrium seeking in networked games," submitted to \emph{Automatica,} revised, available online: arXiv:1912.00415.
\bibitem{HuangCyber18} J. Huang, Y. Song, W. Wang, C. Wen and G. Li, ``Fully distributed adaptive consensus control of a class of high-order nonlinear systems with a directed topology and unknown control directions," \emph{IEEE Transactions on Cybernetics},  vol. 48, no. 8, pp. 2349-2356, 2018.
\bibitem{LiTF14}Z. Li, Z. Duan, \emph{Cooperative Control of Multi-agent Systems: A Consensus Region Approach,} Taylor and Francis/CRC Press, Boca Roton, FL, 2014. ISBN: 978-1-4665-6994-2.

\bibitem{YEAT16} M. Ye, G. Hu, ``A robust extremum seeking scheme for dynamic systems with uncertainties and disturbances," \emph{Automatica}, vol. 66, pp. 172-178, 2016.

\end{thebibliography}
\end{document}